\input amstex
\documentstyle {amsppt}
\UseAMSsymbols \vsize 18cm \widestnumber\key{ZZZZZ}

\catcode`\@=11
\def\displaylinesno #1{\displ@y\halign{
\hbox to\displaywidth{$\@lign\hfil\displaystyle##\hfil$}&
\llap{$##$}\crcr#1\crcr}}
\def\ldisplaylinesno #1{\displ@y\halign{
\hbox to\displaywidth{$\@lign\hfil\displaystyle##\hfil$}&
\kern-\displaywidth\rlap{$##$} \tabskip\displaywidth\crcr#1\crcr}}
\catcode`\@=12

\refstyle{A}

\let \ol=\overline

\let \ti=\widetilde

\font\bf=ecbx1000

\font\nr=eufb7 at 10pt

 \font\srm=cmr10 at 7.5pt

\font\main=cmsy10 at 10pt

\font\smain=cmsy10 at 7.5pt

\font \fin=lasy8 at 15.4 pt
\def \X{\mathop{\hbox{\nr X}^{\hbox{\srm nr}}}\nolimits}

\def \Xim{\mathop{\hbox{\nr X}^{\hbox{\srm nr}}_0}\nolimits}

\def \o{\mathop{\hbox{\main O}}\nolimits}
\def \so{\mathop{\hbox{\smain O}}\nolimits}

\def \F{\mathop{\hbox{\main F}}\nolimits}

\def \End{\mathop{\hbox{\rm End}}\nolimits}

\def \Hom{\mathop{\hbox{\rm Hom}}\nolimits}

\def \deg{\mathop{\hbox{\rm deg}}\nolimits}

\def \det{\mathop{\hbox{\rm det}}\nolimits}

\topmatter
\title Sur quelques r\'esultats de Harish-Chandra pour les repr\'esentations des groupes $p$-adiques, \'elargis aux extensions centrales\endtitle
\rightheadtext{Sur quelques r\'esultats de Harish-Chandra}
\author Volker Heiermann\endauthor
\address Aix-Marseille Universit\'e, CNRS, I2M, 13453 Marseille, France; \endaddress
\email volker.heiermann\@univ-amu.fr \endemail
\abstract
L'objet de cet article est de donner une preuve compl\`ete de r\'esultats de Harish-Chandra liant l'irr\'eductibilit\'e de l'induite parabolique d'une repr\'esentation cuspidale d'un groupe $p$-adique aux propri\'et\'es analytiques de la fonction $\mu $ de Harish-Chandra et de montrer que celle-ci reste valable dans le cas d'une extension centrale.

\endabstract

\endtopmatter
\document

\null Les r\'esultats de Harish-Chandra qui mettent en lien l'irr\'eductibilit\'e de l'induite parabolique d'une repr\'esentation cuspidale d'un groupe $p$-adique et les propri\'et\'es analytiques de la fonction $\mu $ de Harish-Chandra jouent un r\^ole central dans la th\'eorie des repr\'esentations d'un groupe $p$-adique ainsi que dans leur lien avec celle des alg\`ebres de Hecke affines (cf. \cite{H2}, \cite{H3}). Une preuve de ces r\'esultats figure dans l'article \cite{Si2} qui est bas\'e sur le livre \cite{Si1}. Un ingr\'edient important est  la simplicit\'e des p\^oles de la fonction $\mu $ de Harish-Chandra dans le cas cuspidal qui figure dans \cite{Si1}, mais dont la preuve est malheureusement erronn\'ee comme l'auteur l'avait reconnu lui-m\^eme. Ceci avait d\'ej\`a \'et\'e signal\'e dans \cite{H2} lors d'une remarque \`a la preuve de la proposition {\bf 4.1}, tout en mentionnant l'existence d'une autre preuve communiqu\'ee par J.-L. Waldspurger.

L'objet de cette note est de rassembler \`a un seul endroit les preuves de ces r\'esultats de Harish-Chandra y compris celle de la simplicit\'e des p\^oles. La pr\'esenta-tion g\'en\'erale reprend par ailleurs le langage plus r\'ecent de \cite{W}, bien que les preuves soient essentiellement celles dans \cite{Si2}. L'autre motivation est de montrer que tout cela reste valable au cas d'une extension centrale d'un groupe r\'eductif $p$-adique.

\null {\bf 1.} On fixe un corps local non-archim\'edien $F$, un groupe r\'eductif connexe d\'efini sur $F$ dont le groupe des points sera not\'e $G$ et une extension centrale $\ti{G}$ de $G$ donn\'ee par un groupe cyclique d'ordre fini $m$, not\'e $\mu _m$. Nous \'ecrirons $q$ pour le cardinal du corps r\'esiduel de $F$ et $\vert\cdot\vert_F$ pour sa valeur absolue normalis\'ee.

On a donc une suite exacte
$$0\rightarrow{\loadbold\mu_m}\rightarrow^{\iota}\ti{G}\rightarrow^p G\rightarrow 0,$$
l'image de $\mu _m$ \'etant contenue dans le centre de $\ti{G}$. Si $H$ est un sous-ensemble de $G$, on notera $\ti{H}$ son image r\'eciproque dans $\ti{G}$.

On fixe un sous-groupe parabolique minimal $P_0=M_0U_0$ et un sous-groupe parabolique maximal $P=MU$ de $G$, $P\supseteq P_0$. On notera $A_0$ (resp. $A$) le tore d\'eploy\'e maximal dans le centre de $M_0$ (resp. $M$). On fixe par ailleurs un sous-groupe compact maximal $K$ en bonne position relative \`a $A_0$ par la th\'eorie de Bruhat-Tits. On notera $\ol{P}$ le sous-groupe parabolique oppos\'e de $P$, on \'ecrira souvent $\ol{\ti{P}}$ \`a la place de $\ti{\ol{P}}$, et on d\'esignera par $w$ l'unique \'el\'ement du groupe de Weyl de $G$ (relatif \`a $A_0$) tel que $w\ol{P}\supseteq P_0$. (Ici et dans la suite, $w\ol{P}$ d\'esigne l'action du groupe de Weyl sur l'ensemble des sous-groupes de $G$ par conjugaison. Ce sera \'egalement valable pour les sous-groupes de $\ti{G}$.) On a $w\ol{P}=P$, si et seulement si $wM=M$. On fixe par ailleurs un repr\'esentant $\ti{w}$ de $w$ dans $\ti{K}$ tel que $\ti{w}^2\in\ti{M}$. Remarquons que la conjugaison par $\ti{w}$ dans $\ti{G}$ ne d\'epend que de $w$. On a donc de m\^eme $w\ti{M}=\ti{M}$, si et seulement si $w\ti{P}=\ol{\ti{P}}$.

Au moyen d'une section $\ti{P_0}$-\'equivariante $U_0\rightarrow\ti{U_0}$ (cf. \cite{Li2, 2.2}), on a les d\'ecompositions de Levi $\ti{P_0}=\ti{M}U_0$ et $\ti{P}=\ti{M}U$. Cela d\'efinit \'egalement une section pour  $w\ol{U}$, et on prend pour $\ol{U}$ la section conjugu\'ee par $w$. Les foncteurs d'induction parabolique unitaire $i_{\ti{P}}^{\ti{G}}$ et $i_{\ol{\ti{P}}}^{\ti{G}}$ peuvent donc \^etre d\'efinis comme dans le cas d'un groupe r\'eductif, et on notera $r_{\ti{P}}^{\ti{G}}$ et $r_{\ol{\ti{P}}}^{\ti{G}}$ les foncteurs adjoints \`a gauche. Si $(\ti{\pi },\ti{V})$ est une repr\'esentation lisse de $\ti{G}$, on \'ecrira parfois plus simplement $(\ti{\pi}_{\ti{P}},\ti{V}_{\ti{P}})$ et $(\ti{\pi}_{\ol{\ti{P}}},\ti{V}_{\ol{\ti{P}}})$ \`a la place de $(r_{\ti{P}}^{\ti{G}}\ti{\pi},r_{\ti{P}}^{\ti{G}}\ti{V})$ et $(r_{\ol{\ti{P}}}^{\ti{G}}\ti{\pi},r_{\ol{\ti{P}}}^{\ti{G}}\ti{V})$.

On fixe un sous-groupe ouvert ferm\'e  $A^\dagger $ de $A$ v\'erifiant les propri\'et\'es de \cite{Li1, 2.1.1}.  En particulier,  $\ti{A^\dagger}$ est central dans $\ti{A}$ et $\ti{M}$. (On peut par exemple choisir $A^\dagger =A^m$.)

Lorsque $H$ est un sous-groupe alg\'ebrique (resp. un tore) de $G$, $X^*(H)$ (resp. $X_*(H)$) d\'esignera le groupe des caract\`eres (resp. cocaract\`eres) $F$-rationnels de $H$.

On pose $a_{M_0}^*=X^*(M_0)\otimes \Bbb R$ et analogue pour $M$ \`a la place de $M_0$. On v\'erifie ais\'ement que $a_M^*$ s'identifie \`a $X^*(A)\otimes\Bbb R$ par l'homomorphisme de restriction $X^*(M)\rightarrow X^*(A)$. Alors, $a_{M}=X_*(A)\otimes \Bbb R$ est l'espace dual, et on a une d\'ecomposition $a_{M_0}^*=a_M^*\oplus a_{M_0}^{M*},$ o\`u $a_{M_0}^{M*}$ est le sous-espace orthogonal \`a $a_M$. L'espace $a_M^{G*}$ est  donc de dimension $1$, engendr\'e par la restriction \`a $A$ de l'unique racine simple positive de $A_0$ dans $G$ relative \`a $P_0$ qui n'est pas une racine pour $M$. On notera cette racine (ou plut\^ot sa restriction \`a $A$) par $\alpha $ dans la suite. On dira qu'un \'el\'ement de $a_M^{G*}$ est $>0$, s'il est un multiple de $\alpha $ par un nombre r\'eel $>0$.

Pour $\lambda\in\Bbb C$, on notera $\chi _{\lambda\alpha}$ le caract\`ere $m\mapsto \vert \alpha(m)\vert _F^\lambda $ de $M$ et $\ti{\chi}_{\lambda\alpha}$ celui de $\ti{M}$ \'egal \`a $\chi _{\lambda\alpha}\circ p$. Si $\tau $ (resp. $\ti {\tau}$) est une repr\'esentation de $M$ (resp. $\ti{M}$), on \'ecrira $\tau _{\lambda\alpha}$ (resp. $\ti {\tau}_{\lambda\alpha}$) pour la repr\'esentation de $M$ (resp. $\ti{M}$) \'egale \`a $\tau\otimes\chi _{\lambda\alpha}$ (resp.  $\ti {\tau}\otimes\ti{\chi}_{\lambda\alpha}$).

La fonction $\mu $ de Harish-Chandra, $\lambda\mapsto\mu(\ti{\tau }_{\lambda\alpha})$, est d\'efinie dans \cite{W, V.2} pour les groupes $p$-adiques et dans \cite{Li1, 2.4 (7)} pour leurs extensions centrales, de m\^eme l'op\'erateur d'entrelacement $J_{\ol{\ti{P}}\vert \ti{P}}(\ti{\tau}_{\lambda\alpha} ):i_{\ti{P}}^{\ti{G}}\ti{\tau}_{\lambda\alpha} \rightarrow i_{\ol{\ti{P}}}^{\ti{G}}\ti{\tau}_{\lambda\alpha} $. Les deux d\'ependent m\'eromorphiquement de $\lambda $. En fait, ils peuvent \^etre vus comme applications rationnelles comme on le rappelera ci-dessous.

Le r\'esultat suivant est d\^u \`a Harish-Chandra \cite{Si1, Si2} dans le cas des groupes $p$-adiques. (En fait, comme d\'ej\`a rappel\'e dans l'introduction, seule la simplicit\'e du p\^ole est \`a reprendre.)

\null{\bf 1.1 Th\'eor\`eme:} \it Soit $\ti{\sigma} $ une repr\'esentation irr\'eductible cuspidale unitaire de $\ti{M}$.

Alors,

i) Si $w\ti{M}\ne \ti{M}$ ou $w\ti{\sigma}\not\simeq\ti{\sigma}$, alors $i_{\ti{P}}^{\ti{G}}\ti{\sigma}_{\lambda\alpha }$ est irr\'eductible pour tout $\lambda\in\Bbb R$ et l'application $\lambda\mapsto\mu(\ti{\sigma}_{\lambda\alpha })$ est r\'eguli\`ere et non nulle pour $\lambda\in\Bbb R$.

ii) Supposons $w\ti{M}=\ti{M}$ et $w\ti{\sigma}\simeq\ti{\sigma} $. Alors,

a) il existe un unique $\lambda\geq 0$ tel que $i_{\ti{P}}^{\ti{G}}\ti{\sigma}_{\lambda\alpha}$ soit r\'eductible, et alors $i_{\ti{P}}^{\ti{G}}\ti{\sigma}_{-\lambda\alpha}$ est \'egalement r\'eductible.

b) la repr\'esentation $i_{\ti{P}}^{\ti{G}}\ti{\sigma} $ est irr\'eductible, si et seulement si $\mu(\ti{\sigma} )=0$. Si $\mu(\ti{\sigma} )\ne 0$, elle est la somme de deux repr\'esentations irr\'eductibles non isomorphes.

c) soit $\lambda >0$. Alors, la fonction inverse de $\mu $ est d\'efinie en $\ti{\sigma}_{\lambda\alpha}$. La repr\'esenta-tion $i_{\ti{P}}^{\ti{G}}\ti{\sigma}_{\lambda\alpha }$ est r\'eductible, si et seulement si $\lambda$ est un p\^ole de l'application $\lambda\mapsto\mu(\ti{\sigma} _{\lambda\alpha })$. Ce p\^ole est simple.

(iii) Les p\^oles de l'op\'erateur d'entrelacement $J_{\ol{\ti{P}}\vert \ti{P}}$ sur l'ensemble des $\ti{\sigma}_{\lambda\alpha}$, $\lambda\in\Bbb C$, sont exactement les z\'eros de la fonction $\mu $. Ils sont tous simples.

(iv) Si l'op\'erateur d'entrelacement $J_{\ol{\ti{P}}\vert \ti{P}}$ est r\'egulier en $\ti{\sigma}$, alors il est r\'egulier et bijectif en $\ti{\sigma}_{\lambda\alpha}$ pour tout $\lambda\in\Bbb R$.

\rm

\null
Notons $M^1=\cap_{\chi \in X^*(M)} \ker(\vert \chi\vert_F )$, $G^1=\cap_{\chi \in X^*(G)} \ker(\vert \chi\vert_F )$ et posons $B=\Bbb C[M/M^1]$. \'Ecrivons $\X(M)$ pour le groupe des homomorphismes $M\rightarrow\Bbb C^*$ dont le noyau contient $M^1$, $\X(G)$ pour le groupe des homomorphismes $G\rightarrow\Bbb C^*$ dont le noyau contient $G^1$, $\X(\ti{M})=\X(M)\circ p$ et $\X(\ti{G})=\X(G)\circ p$. Les caract\`eres $\ti{\chi }_{\lambda \alpha}$, $\lambda\in\Bbb C$, sont triviaux sur le centre de $\ti{G}$ et appartiennent \`a $\X (\ti{M})$. L'anneau $B$ s'identifie \`a l'anneau des polyn\^omes sur la vari\'et\'e $\X (\ti{M})$. Pour $m\in M$, $b_m$ d\'esignera l'\'el\'ement de $B$ support\'e en $mM^1$. Pour $m\in\ti{M}$, on \'ecrira par abus de notation $b_m=b_{p(m)}$. On d\'efinit $b_m(\ti{\chi })=\ti{\chi }(m)$. Il a \'et\'e prouv\'e dans \cite{W, IV.1.1} pour les groupes $p$-adiques (et v\'erifi\'e dans \cite{Li1, 2.4.1} pour leurs extensions centrales) que la fonction $\mu(\ti{\sigma}\otimes\cdot)$ et l'op\'erateur d'entrelacement $J_{\ol{\ti{P}}\vert \ti{P}}(\ti{\sigma}\otimes\cdot)$ sont rationnels en $\ti{\chi}\in\X(\ti{M})$. En fait, ils ne d\'ependent que des classes modulo\footnote{Cela se voit sur l'expression de l'op\'erateur d'entrelacement par une int\'egrale (cf. {\bf 2.2}) \`a l'aide de la d\'ecomposition d'Iwasawa et du fait que les caract\`eres de $\ti{G}$ sont triviaux sur le groupe d\'eriv\'e qui contient les \'el\'ements unipotents.} $\X(\ti{G})_{\vert\ti{M}}$ et en un certain sense que de la classe d'isomorphie de $\ti{\sigma }$ (cf. \cite{W, IV.1}).

On a choisi dans \cite{H3, 1.5} un \'el\'ement particulier $h_{\alpha}\in M$ associ\'e \`a $\alpha $ et d\'efini $Y_{\alpha }=b_{h_{\alpha}}\in B$. On a $\ti{\chi}(h_{\alpha})=1$, si et seulement si $\ti{\chi}\in\X(\ti{G})$. Il existe un plus petit entier $t_{\alpha }>0$ tel que, pour $\ti{\chi}\in \X(\ti{M})$,  $\ti{\chi }(h_{\alpha }^{t_{\alpha }})=1$ \'equivaut \`a l'existence de $\ti{\chi'}\in\X(\ti{G})$ v\'erifiant $\ti{\sigma}\otimes\ti{\chi}\ti{\chi'}\simeq\ti{\sigma}$. Posons $X_{\alpha }=Y_{\alpha }^{t_{\alpha }}$. On \'ecrira ici simplement $Y$ et $X$. Remarquons que $X(\ti{\chi} _{\lambda\alpha})=q^{-t\lambda}$ pour un certain nombre r\'eel $t>0$ et que la valeur $X(\ti{\chi })$ d\'etermine par d\'efinition la classe d'isomorphie de $\ti{\sigma }\otimes\ti{\chi }$ \`a torsion par un caract\`ere de $\X(\ti{G})$ pr\`es.

\null{\bf 1.2 Lemme:} \it Avec les notations du th\'eor\`eme, supposons $w\ti{M}=\ti{M}$ et $w\ti{\sigma}\simeq\ti{\sigma} $. Si $\ti{\chi}\in\X(\ti{M})$ v\'erifie $w(\ti{\sigma }\otimes\ti{\chi})\simeq\ti{\sigma }\otimes\ti{\chi }$, alors $X(\ti{\chi })\in\{\pm 1\}$. R\'eciproquement, l'isomorphie est toujours v\'erifi\'e si $X(\ti{\chi })=1$. En outre, \`a isomorphisme et torsion avec un caract\`ere de $\X(\ti{G})$ pr\`es, il existe au plus une repr\'esentation irr\'eductible cuspidale $\ti{\sigma}_-$ de la forme $\ti{\sigma }\otimes\chi_{\lambda\alpha }$, $\lambda\in\Bbb C$, tel que $X(\chi_{\lambda\alpha })=-1$ et $w\ti{\sigma}_-\simeq \ti{\sigma}_-$. Elle est unitaire \`a torsion avec un caract\`ere de $\X(\ti{G})$ pr\`es.

\null Preuve: \rm Si $X(\ti{\chi })=1$, alors il existe par d\'efinition $\ti{\chi'}\in\X(\ti{G})$ tel que $\ti{\sigma }\otimes\ti{\chi} {\ti{\chi'}}^{-1}\simeq\ti{\sigma }$ et donc
$$w(\ti{\sigma }\otimes\ti{\chi})=w(\ti{\sigma }\otimes\ti{\chi{\chi'}}^{-1}\ti{\chi'})\simeq w(\ti{\sigma }\otimes\ti{\chi'})=(w\ti{\sigma })\otimes\ti{\chi'}\simeq \ti{\sigma }\otimes\ti{\chi'}\simeq\ti{\sigma }\otimes\ti{\chi}\ti{{\chi'}^{-1}}\ti{\chi'}=\ti{\sigma }\otimes\ti{\chi}.$$

Supposons maintenant $X(\ti{\chi })\ne 1$ et \'ecrivons $\ti{\chi}=\ti{\chi'}\ti{\chi}_{\lambda\alpha}$ avec $\ti{\chi '}\in\X(\ti{G})$ et $\lambda\in\Bbb C$. Alors, $$w(\ti{\sigma}\otimes\ti{\chi'}\ti{\chi}_{\lambda\alpha})\simeq\ti{\sigma}\otimes\ti{\chi'}\ti{\chi}_{-\lambda\alpha},$$
et, par cons\'equence, $w(\ti{\sigma}\otimes\ti{\chi})\simeq\ti{\sigma}\otimes\ti{\chi}$ \'equivaut \`a  $\ti{\sigma}\otimes\ti{\chi}_{\lambda\alpha}\simeq \ti{\sigma}\otimes\ti{\chi}_{-\lambda\alpha}$ ou encore \`a $\ti{\sigma}\otimes\ti{\chi}_{2\lambda\alpha}\simeq \ti{\sigma}$, d'o\`u $1=X(\ti{\chi}_{2\lambda\alpha})=X(\ti{\chi}_{\lambda\alpha})^2$ et $\lambda\in i\Bbb R$. En particulier, $\ti{\sigma }\otimes\ti{\chi}$ est unitaire \`a torsion avec un caract\`ere de $\X(\ti{G})_{\vert \ti{M}}$ pr\`es. Comme par hypoth\`ese $X(\ti{\chi}_{\lambda\alpha})=X(\ti{\chi })\ne 1$, on a bien $-1=X(\ti{\chi}_{\lambda\alpha})=X(\ti{\chi })$. Cela d\'etermine $\ti{\sigma }\otimes\ti{\chi}$ \`a torsion avec un caract\`ere de $\X(\ti{G})$ pr\`es, comme d\'ej\`a remarqu\'e au-dessus de l'\'enonc\'e. \hfill{\fin 2}

\null {\bf 1.3} \it Remarque: \rm  Ayant $X(\ti{\chi} _{\lambda\alpha})=q^{-t\lambda}$ avec $t$ un nombre r\'eel  $>0$, on peut toujours trouver $\lambda_0\in\Bbb C$ tel que $X(\ti{\chi }_{\lambda _0\alpha})=-1$, mais cela n'implique pas n\'ecessairement $w(\ti{\sigma}\otimes\ti{\chi})\simeq\ti{\sigma}\otimes\ti{\chi}$ et alors aucun $\ti{\chi }$ avec $X(\ti{\chi })=-1$ n'a cette propri\'et\'e: prenons par exemple $G=GL_{2n}(F)$, $M=GL_n(F)\times GL_n(F)$ et $\ti{\sigma }=\sigma_0\otimes\sigma_0$, $\sigma _0$ repr\'esentation cuspidale unitaire de $GL_n(F)$. Alors, $\sigma_0\vert\det_n\vert ^x\otimes\sigma_0\vert\det_n\vert ^{-x}$ est invariant par $w$ si et seulement si $\sigma_0\vert\det_n\vert ^{2x}\simeq\sigma_0$, mais alors $\sigma_0\vert\det_n\vert ^x\otimes\sigma_0\vert\det_n\vert ^{-x}\simeq \sigma_0\vert\det_n\vert ^x\otimes\sigma_0\vert\det_n\vert ^{x}\simeq (\sigma_0\otimes\sigma_0) \vert\det_{2n}\vert ^x$. Donc, $X(\vert\det_n \vert ^x\otimes\vert\det_n\vert ^{-x})=1$ par d\'efinition de $X$.

\null Il est clair par la rationnalit\'e remarqu\'ee ci-dessus que $\mu (\ti{\sigma}\otimes\cdot)$ est rationnelle en $X$ et que $J_{\ol{\ti{P}}\vert \ti{P}}(\ti{\sigma}\otimes\cdot)$ est rationnelle en $Y$. On a un r\'esultat plus pr\'ecis qui r\'esulte du th\'eor\`eme {\bf 1.1} (pour les groupes $p$-adiques, c'est dans \cite{Si2}):

\null{\bf 1.4 Corollaire:}  \it Soit $\ti{\sigma}$ une repr\'esentation cuspidale unitaire de $\ti{M}$ et supposons que $w\ti{\sigma}\simeq\ti{\sigma}$ s'il existe $\ti{\chi}\in\X(\ti{M})$ tel que $\ti{\sigma}\otimes\ti{\chi}$ v\'erifie cette propri\'et\'e. Alors, il y a des nombres r\'eels $a, a_-\geq 0$ et $c>0$ tels que
$$\mu (\ti{\sigma}\otimes\cdot)=c\ {(1-X(\cdot ))(1-X^{-1})(\cdot )\over (1-X(\cdot )q^{-a})(1-X^{-1}(\cdot )q^{-a})} {(1+X(\cdot ))(1+X^{-1})(\cdot )\over (1+X(\cdot )q^{-a_-})(1+X^{-1}(\cdot )q^{-a_-})}.$$

\null Preuve: \rm On sait que $\mu $ est une fonction rationnelle en $X$. Choisissons une repr\'esentation cuspidale unitaire $\ti{\sigma}_-$ conforme au lemme {\bf 1.2} (s'il y en a). On voit que $a$ correspond \`a l'unique p\^ole $>0$ (s'il y en a) sur la droite r\'eelle $\ti{\sigma}_{\lambda\alpha}$ et $b$ correspond \`a l'unique p\^ole $>0$ (s'il y en a) sur la droite r\'eelle $(\ti{\sigma}_-)_{\lambda\alpha}$. Par le th\'eor\`eme {\bf 1.1}, le lemme {\bf 1.2} et l'invariance de $\mu (\ti{\sigma}\otimes\cdot)$ par torsion par $\X(\ti{G})$, ils ne peuvent pas y avoir d'autres p\^oles.

S'il n'y a pas de p\^oles, alors, par le th\'eor\`eme {\bf 1.1}, il n'y a pas non plus de z\'eros, i.e. $\mu $ est un mon\^ome en $X$. Mais, comme les valeurs de $\mu $ sur l'axe unitaire sont r\'eelles et $\geq 0$ (cf. \cite{W, V.2.1}, \cite{Li1, 2.4.3}) et que $\mu $ n'est pas non plus identiquement nulle, on trouve bien que $\mu \equiv c$ avec $c$ r\'eel $>0$ et on a $a=a_-=0$.

De m\^eme, le th\'eor\`eme {\bf 1.1} montre que, s'il n'y a pas de p\^ole  sur la droite r\'eelle $\ti{\sigma}_{\lambda\alpha}$ (resp. sur la droite r\'eelle $(\ti{\sigma}_-)_{\lambda\alpha}$), il n'y a pas non plus de z\'ero sur cette droite, et inversement. On pose $a=0$ (resp. $a_-=0$), s'il n'y a pas de z\'ero.

Ceci prouve, compte tenu de la multiplicit\'e des p\^oles et des z\'eros, que la fonction $\mu $ a la forme indiqu\'ee \`a multiplication avec un mon\^ome pr\`es avec un coefficient complexe de valeur absolue $1$. Toutefois, l'expression \`a droite dans l'\'egalit\'e de l'\'enonc\'e est $\geq 0$ et non identiquement nulle sur l'axe unitaire comme l'est la fonction $\mu $. Mais, pour qu'un mon\^ome de la forme $c'X^k$ avec $c'$ complexe de valeur absolue $1$ et $k\in\Bbb Z$ prenne uniquement des valeurs r\'eelles $>0$ sur l'axe unitaire, il faut et il suffit que $c'=1$ et $k=0$.\hfill{\fin 2}

\null La preuve du th\'eor\`eme {\bf 1.1} se fera dans la section {\bf 4.} apr\`es plusieurs r\'esultats interm\'ediaires.

\null{\bf 2.} Fixons une repr\'esentation irr\'eductible cuspidale unitaire $(\ti{\sigma },\ti{E})$ de $\ti{M}$ et notons $\chi_{\ti{\sigma }}$ son caract\`ere central.  Suivant \cite{W, I.5}, posons $\ti{E}_B=\ti{E}\otimes B$ et notons $\ti{\sigma }_B$ la repr\'esentation de $\ti{M}$ dans $\ti{E}_B$ donn\'ee par $\ti{\sigma }_B(m)(e\otimes b)=\ti{\sigma }(m)e\otimes bb_m$. Remarquons que  $\ti{E}_B$ est un $(\ti{M},B)$-module, l'action de $\ti{M}$ et l'action naturelle de $B$ sur le deuxi\`eme facteur commutant.

On utilisera dans {\bf 2.5, 2.6} et dans les annexes que les repr\'esentations induites $i_{\ti{P}}^{\ti{G}}\ti{\sigma}_{\lambda\alpha }$, $\lambda\in\Bbb C$, puissent toutes \^etre r\'ealis\'ees dans le m\^eme espace $i_{\ti{P}\cap \ti{K}}^{\ti{K}}\ti{E}$.

\null {\bf 2.1 Proposition:} \it La semisimplifi\'ee de $r_{\ti{P}}^{\ti{G}}(i_{\ti{P}}^{\ti{G}}\ti{\sigma}_{\lambda\alpha })$ (resp. $r_{\ol{\ti{P}}}^{\ti{G}}(i_{\ti{P}}^{\ti{G}}\ti{\sigma}_{\lambda\alpha })$) est $\ti{\sigma}_{\lambda\alpha  }$ si $w\ti{M}\ne \ti{M}$, et $\ti{\sigma}_{\lambda\alpha }\oplus (w\ti{\sigma})_{-\lambda\alpha }$ si $w\ti{M}=\ti{M}$.

\null Preuve: \rm Cela r\'esulte de \cite{BZ, 5.2}.\hfill{\fin 2}

\null{\bf 2.2} L'op\'erateur d'entrelacement $J_{\ol{\ti{P}}\vert \ti{P}}(\ti{\sigma}_{\lambda \alpha})$ est d\'efini pour $\lambda >0$ par  $$i_{\ti{P}}^{\ti{G}}\ti{\sigma}_{\lambda \alpha}\rightarrow i_{\ol{\ti{P}}}^{\ti{G}}\ti{\sigma}_{\lambda \alpha}, v\mapsto \{g\mapsto \int_{\ol{U}} v(ug)du\}$$ et pour tout $\lambda $ par prolongement analytique, remarquant que les coefficients matriciels sont des fonctions rationnelles en $Y$. La preuve dans \cite{W} se g\'en\'eralise au cas des extensions centrales avec les notations dans \cite{Li1}, rempla\c cant $A$ par $\ti{A^\dagger}$. Lors de cette preuve  \cite{W,  IV.1.1, p. 282-3}, il est notamment prouv\'e le r\'esultat suivant:

\null {\bf Lemme:} \it Supposons $w\ti{P}=\ol{\ti{P}}$. Posons $\F _1=\{v\in i_{\ti{P}}^{\ti{G}}\ti{E}_B\vert supp(v)\subseteq \ti{P}\ol{\ti{P}}\}$. Fixons $a\in\ti{A^\dagger}$ tel que $waw^{-1}\ne a$ et notons $b=\chi_{\ti{\sigma} }(a)b_a- (w\chi _{\ti{\sigma} })(a)\ ^wb_a$.

L'endomorphisme d'espaces vectoriels $p_1:\F_1\rightarrow \ti{E}_B$, $v\mapsto (J_{\ol{\ti{P}}\vert\ti{P}}(\ti{\sigma} _{B})v)(1)$ est bien d\'efini. Il se factorise en un homomorphisme de $(\ti{M},B)$-modules $q_1:(\F_1)_{\ol{\ti{P}}}\rightarrow \ti{E}_B$, $\ol{v}\mapsto (J_{\ol{\ti{P}}\vert\ti{P}}(\ti{\sigma} _B)v)(1)$.

L'op\'erateur $bJ_{\ol{\ti{P}}\vert\ti{P}}$ appartient \`a $Hom_{\ti{G},B}(i_{\ti{P}}^{\ti{G}}\ti{E}_B, i_{\ol{\ti{P}}}^{\ti{G}}\ti{E}_B)$ et vient par r\'eciprocit\'e de Frobenius de l'\'el\'ement $j\in Hom_{\ti{G},B}((i_{\ti{P}}^{\ti{G}}\ti{E}_B)_{\ol{\ti{P}}},\ti{E}_B)$ \'egal \`a $q_1\circ (r_{\ol{\ti{P}}}^{\ti{G}}(i_{\ti{P}}^{\ti{G}}\ti{\sigma} _B)(a-(w\chi_{\ti{\sigma} })(a)\ ^wb_a))$.\rm

\null Remarquons que $b_a=Y^{t_a}$ pour un certain entier $t_a$ et que l'on peut choisir $a$ tel que $t_a>0$. Alors, le $b$ dans l'\'enonc\'e du lemme s'\'ecrit $\chi_{\ti{\sigma} }(a)Y^{t_a}- (w\chi _{\ti{\sigma} })(a) Y^{-t_a}$.

\null {\bf 2.3 Corollaire:} \it L'op\'erateur d'entrelacement $J_{\ol{\ti{P}}\vert\ti{P}}$ est r\'egulier non nul en tout $\ti{\sigma }_{\lambda\alpha }$ avec $\lambda $ r\'eel non nul. Il est r\'egulier en $\ti{\sigma} $ si $w\ti{M}\ne \ti{M}$ ou si $w\ti{\sigma} \not\simeq\ti{\sigma} $. Si $w\ti{M}=\ti{M}$ et $w\ti{\sigma}\simeq\ti{\sigma}$, alors cet op\'erateur a au plus un p\^ole simple en $\ti{\sigma}$.

\null Preuve: \rm Pour les deux premiers points, on renvoie \`a \cite{W, IV.2.1}, \cite{W, IV.1.2} et \cite{W, IV.1.1 (10)} qui sont bas\'es sur le lemme pr\'ec\'edent ou son analogue si $w\ti{P}\ne\ol{\ti{P}}$ et se g\'en\'eralisent sans probl\`eme aux extensions centrales. Le deuxi\`eme point r\'esulte directement du lemme. \hfill{\fin 2}

\null{\bf 2.4} Remarquons que la fonction $\mu $ de Harish-Chadra est ici donn\'ee par la formule $\gamma(\ti{G}/\ti{M})^2\mu(\ti{\sigma }_{\lambda\alpha })^{-1}=J_{\ti{P}\vert \ol{\ti{P}}}(\ti{\sigma }_{\lambda\alpha })J_{\ol{\ti{P}}\vert \ti{P}}(\ti{\sigma }_{\lambda\alpha })$ avec $\gamma(\ti{G}/\ti{M})$ une constante (cf. \cite{W, V.2, I.1 (3)}.\footnote{Cette constante ne figure pas dans \cite{Li1} (cf. \cite{Li1, 2.4 (7)}) qui a choisi une autre convention pour les constantes dans la formule de Plancherel. Elle avait \'egalement d\'ej\`a \'et\'e omise dans \cite{H2} et \cite{H3}. Les notations dans \cite{W} ont \'et\'e bas\'ees sur un manuscrit de Harish-Chandra sur la formule de Plancherel pour les groupes $p$-adiques.}

Ceci fait bien du sens: le produit des deux op\'erateurs d'entrelacement est certainement scalaire par le lemme de Schur si la repr\'esentation induite est irr\'eductible. Par ailleurs, la compos\'ee est une fonction rationnelle sur $\X(\ti{M})$. Choisissons un sous-groupe ouvert compact $\ti{H}$ comme dans \cite{BD, 3.9 (ii)} qui s'applique \'egalement aux extensions centrales (cf. remarque en bas de la page 16 de \cite{BD}). L'irr\'eductibi-lit\'e peut alors se lire sur le sous-espace des $\ti{H}$-invariants qui est de dimension finie. L'op\'erateur d'entrelacement y correspond \`a une matrice qui est, par ce qui pr\'ec\`ede, scalaire si la repr\'esentation est irr\'eductible. Par les lemmes {\bf 3.1} et {\bf 3.2}  ci-dessous (parties qui n'utilisent pas la fonction $\mu $), la repr\'esentation induite est irr\'eductible dans l'ouvert de Zarisky donn\'e par la non-nullit\'e du num\'erateur du d\'eterminant de cette matrice. Par cons\'equence, la compos\'ee est scalaire sur un ouvert de Zarisky, donc scalaire partout par prolongement analytique. Ceci \'etant vrai pour des $\ti{H}$ arbitrairement petits , la compos\'ee est bien un op\'erateur scalaire.

\null \null{\bf Proposition:} \it Supposons $\lambda$ non nul et que $\ti{\sigma} _{\lambda\alpha }$ soit un p\^ole de $\mu $.

(i) La repr\'esentation induite $i_{\ti{P}}^{\ti{G}}\ti{\sigma} _{\lambda\alpha }$ est alors r\'eductible non semi-simple de longueur $2$ et poss\`ede un unique sous-quotient de carr\'e int\'egrable. C'est l'unique sous-repr\'esentation irr\'eductible si $\lambda>0$.

(ii) Si $\lambda>0$, alors $Res_{z=\lambda}\mu(\ti{\sigma}_{z\alpha })>0$.\rm

\null\it Preuve: \rm Si $\mu $ a un p\^ole en $\ti{\sigma}_{\lambda\alpha}$, $\lambda >0$, alors la compos\'ee des deux op\'erateurs d'entrelacement est nulle. Comme les deux sont r\'eguliers par le corollaire {\bf 2.3}, au moins un d'eux ne peut pas \^etre bijectif. Mais, comme il est non nul par {\bf 2.3}, l'induite admet alors un sous-espace propre non nul invariant. On a donc r\'eductibilit\'e. Il r\'esulte alors de la r\'eciprocit\'e de Frob\'enius avec la proposition {\bf 2.1} que cette induite est de longueur $2$ et non semi-simple, et du crit\`ere de Casselman (cf. \cite{W, III.1.1} et \cite{Li1, 2.3.1}) qu'il poss\`ede un unique sous-quotient irr\'eductible de carr\'e int\'egrable qui est l'unique sous-repr\'esentation irr\'eductible si $\lambda>0$.

Il est montr\'e dans la proposition {\bf A.4} de l'annexe que, si $\lambda>0$, le degr\'e formel de la sous-repr\'esentation irr\'eductible de carr\'e int\'egrable de $\ti{\sigma} _{\lambda\alpha }$ est le multiple par ${\iota_M\over\iota_G}\log q\ \gamma(\ti{G}\vert\ti{M})^{-1} Res_{z=0}\mu(\ti{\sigma}_{(z+\lambda)\ti{\alpha }})$ du degr\'e formel de $\ti{\sigma }$. Le degr\'e formel d'une repr\'esentation de carr\'e int\'egrable \'etant $>0$ (cf. \cite{W, III.1 p. 265} et \cite{Li1, 2.3}) ainsi que les constants ${\iota_G\over\iota_M}$, $\log q$, $\gamma(\ti{G}\vert\ti{M})$ et le ratio $\alpha/\ti{\alpha}$, on en d\'eduit le (ii).\hfill{\fin 2}

\null{\bf 2.5} La proposition ci-dessous sera utilis\'ee seulement dans la preuve du lemme {\bf 3.4} (ii) ainsi que dans celle du th\'eor\`eme {\bf A.2} dans l'annexe. Elle r\'esulte de la formule d'approximation de Casselman \cite{H1, 1.3.1}, bas\'ee sur l'accouplement de Casselman \cite{W, I.4.1}. Cela reste valable pour les extensions centrales si on passe \`a $\ti{A^\dagger}$: en effet, pour l'accouplement de Casselman cela a \'et\'e v\'erifi\'e dans \cite{KS, preuve de 3.10}, alors que les autres arguments sont en bonne partie combinatoires utilisant les propri\'et\'es des op\'erateurs d'entrelacement ou du foncteur de Jacquet qui restent vraies dans le cas d'une extension centrale.

Notons $(\ti{\sigma}^\vee,\ti{E}^\vee)$ la repr\'esentation duale de $(\ti{\sigma},\ti{E})$.

\null{\bf Proposition:} \it Supposons $w\ti{M}=\ti{M}$ et $w\ti{\sigma }=\ti{\sigma }$. Notons $\lambda(\ti{w})$ l'isomorphisme $i_{\ti{P}}^{\ti{G}}\ti{\sigma }\rightarrow i_{\ti{wP}}^{\ti{G}}\ti{w}\ti{\sigma }$ d\'efini par translation \`a gauche par $\ti{w}$ dans l'espace des fonctions $\ti{G}\rightarrow \ti{E}$ dans l'espace induit.

Soient $v\in i_{\ti{P}\cap \ti{K}}^{\ti{K}}\ti{E}$ et $v^{\vee }\in i_{\ti{P}\cap \ti{K}}^{\ti{K}}\ti{E}^{\vee }$. Il existe $\epsilon >0$, de sorte que, pour tout $\chi \in \X (\ti{M})$ et tout $a\in \ti{A^\dagger}$, tel que $\vert \alpha(a)\vert\leq\epsilon $, on ait
$$\eqalign {
&\langle i_{\ti{P}}^{\ti{G}}(\ti{\sigma}_{\lambda\alpha })(a)v,v^{\vee }\rangle\cr
=& \gamma (\ti{G}\vert{\ti{M}})^{-1}\delta _{\ti{P}}^{1/2}(a)\ \{\chi_{\ti{\sigma }}(a)\ \vert\alpha(a)\vert^{\lambda}\ \langle v_{\vert{\ti{M}}}, (J_{\ol{\ti{P}}\vert \ti{P}}(\ti{\sigma}^\vee _{-\lambda\alpha })v^{\vee })_{\vert{\ti{M}}}\rangle_{\ti{M}}\cr &\qquad\qquad\qquad\ \ \ \ \ \ +\chi_{w\ti{\sigma}}(a)\ \vert\alpha(a)\vert^{-\lambda}\ \langle (\lambda(\ti{w})J_{\ol{\ti{P}}\vert \ti{P}}(\ti{\sigma} _{\lambda\alpha })v)_{\vert{\ti{M}}}, v^{\vee }_{\vert{\ti{M}}}\rangle_{\ti{M}}\}.
  \cr}$$
En particulier, si $\lambda\alpha>\delta_{\ti{P}}^{1/2}$, alors $\vert <i_{\ti{P}}^{\ti{G}}(\ti{\sigma}_{\lambda\alpha })(a)v, v^{\vee }> \vert$ tend vers l'infini lorsque $\vert\alpha(a)\vert_F\rightarrow 0$.
\rm

\null {\bf 2.6} Fixons une forme hermitienne $\ti{M}$-invariante d\'efinie positive $\langle\cdot,\cdot\rangle_{\ti{\sigma}}$ sur l'espace de $\ti{\sigma} $. Supposons $w\ti{\sigma}\simeq\ti{\sigma} $. Observons que $i_{\ti{P}}^{\ti{G}}(\ti{\sigma}_{\lambda\alpha })$ est alors hermitienne puisque isomorphe \`a $i_{\ti{P}}^{\ti{G}}(\ti{\sigma}_{-\lambda\alpha })$ qui s'identifie au dual hermitien de $i_{\ti{P}}^{\ti{G}}(\ti{\sigma}_{\lambda\alpha })$ via la forme sesquilin\'eaire $\ti{G}$-\'equivariante $$\langle v,v'\rangle_{\ti{P},\ti{\sigma }}:=\int _{\ti{P}\backslash \ti{G}} \langle v(g), v'(g)\rangle_{\ti{\sigma}} dg$$ qui d\'efinit par ailleurs un produit scalaire sur l'espace vectoriel $i_{\ti{P}\cap\ti{K}}^{\ti{K}}\ti{E}$.

Posons $\epsilon =1$ ou $0$ selon que $J_{\ol{\ti{P}}\vert\ti{P}}$ ait un p\^ole en $\ti{\sigma }$ ou pas. Rappelons que celui-ci est au plus simple (cf. {\bf 2.3}). L' op\'erateur $$(X-1)^\epsilon\lambda(\ti{w})J_{\ol{\ti{P}}\vert\ti{P}}(\ti{\sigma }_{\lambda\alpha })\eqno{\hbox{\rm (*)}}$$ est donc r\'egulier et non nul en tout $\lambda\in\Bbb R$ et d\'efinit un homomorphisme $\ti{G}$-\'equivariant $i_{\ti{P}}^{\ti{G}}\ti{\sigma }_{\lambda\alpha }\rightarrow i_{\ti{P}}^{\ti{G}}(\ti{w}\ti{\sigma })_{-\lambda\alpha }$. On a d\'efini dans l'annexe {\bf B.} suivant \cite{H3, 2.4 et 2.5} un isomorphisme $\rho_{\ti{\sigma },\ti{w}}:\ti{w}\ti{\sigma }\rightarrow \ti{\sigma }$ et prouv\'e que l'op\'erateur $R_{\ti{P}}(\ti{\sigma}_{\lambda\alpha },\ti{w})$ obtenu en composant l'op\'erateur (*) \`a gauche avec $\rho_{\ti{\sigma },\ti{w}}$ est autoadjoint pour $\langle\cdot ,\cdot\rangle_{\ti{P},\ti{\sigma }}$.

\null {\bf Proposition:} \it Supposons $w\ti{\sigma}\simeq\ti{\sigma} $. Soit $\lambda\geq 0$.
La forme sesquilin\'eaire $\ti{G}$-\'equivariante $$(v,v')\mapsto \langle v,(R_{\ti{P}}(\ti{\sigma} _{\lambda\alpha },\ti{w})v'\rangle_{\ti{P},\ti{\sigma }}$$ d\'efinie sur l'espace de $i_{\ti{P}}^{\ti{G}}\ti{\sigma}_{\lambda\alpha }$ induit pour $\lambda >0$, par passage au quotient, un produit hermitien sur le quotient de Langlands de $i_{\ti{P}}^{\ti{G}}\ti{\sigma}_{\lambda\alpha }$. Cette forme est \'egalement hermitienne si $\lambda =0$. Pour $\lambda>0$, le quotient de Langlands est unitaire, si et seulement si ce produit hermitien est d\'efini (positif ou n\'egatif).

\null Preuve: \rm Par la r\'eciprocit\'e de Frobenius et {\bf 2.1}, la repr\'esentation induite $i_{\ti{P}}^{\ti{G}}\ti{\sigma}_{\lambda\alpha }$ est ou bien irr\'eductible ou bien de longueur $2$. Donc, l'op\'erateur d'entrelacement \'etablit, pour $\lambda >0$, un isomorphisme entre le quotient de Langlands\footnote{Ceux qui sont incomfortable avec l'utilisation de la th\'eorie du quotient de Langlands ici peuvent consid\'erer \`a la place la proposition {\bf 3.2} qui ne d\'epend pas de l'\'enonc\'e \`a prouver.} qui est irr\'eductible et le dual hermitien de celui-ci. Ainsi, on obtient bien une forme sesquilin\'eaire $\ti{G}$-invariante sur ce quotient irr\'eductible. Elle est hermitienne, puis-que $R_{\ti{P}}(\ti{\sigma}_{\lambda\alpha },\ti{w})$ est autoadjoint. Une telle forme sur une repr\'esentation irr\'eductible est uniquement d\'etermin\'ee \`a une constante pr\`es. Donc, l'unitarit\'e \'equivaut \`a dire que le produit est d\'efini positif ou n\'egatif.

Lorsque $\lambda =0$, il r\'esulte de ce qui pr\'ec\'edait par prolongement analytique que la forme est toujours hermitienne. \hfill{\fin 2}

\null{\bf 3.} Cette section est consacr\'ee \`a la d\'emonstration de r\'esultats pr\'eliminaires \`a la preuve du th\'eor\`eme {\bf 1.1}.

\null {\bf 3.1 Lemme:} \it La repr\'esentation $i_{\ti{P}}^{\ti{G}}\ti{\sigma} $ est irr\'eductible, si et seulement si $w\ti{M}\ne \ti{M}$, $w\ti{\sigma}\ne \ti{\sigma} $ ou si l'op\'erateur d'entrelacement $J_{\ol{\ti{P}}\vert \ti{P}}$ a un p\^ole en $\ti{\sigma} $. Ce p\^ole est alors simple.

Supposons $i_{\ti{P}}^{\ti{G}}\ti{\sigma} $ r\'eductible. Alors, cette repr\'esentation est somme directe de deux repr\'esentations irr\'eductibles non-isomorphes. Fixons un isomorphisme $\rho _w:w\ti{\sigma }\rightarrow\ti\sigma $ et rappelons l'isomorphisme $\lambda(\ti{w})$ d\'efini dans la proposition {\bf 2.5}. L'op\'erateur $\rho_w\lambda(\ti{w})J_{\ol{\ti{P}}\vert \ti{P}}(\ti{\sigma} )$ est un automorphisme de $i_{\ti{P}}^{\ti{G}}\ti{\sigma} $ qui agit sur chacune des sous-repr\'e-sentations irr\'eductibles par un scalaire, ces deux scalaires diff\'erant par un facteur $-1$.

\null Remarque: \rm On verra lors de la preuve que l'on peut normaliser l'isomorphisme $\rho _w$ de telle sorte que le scalaire soit de carr\'e \'egal \`a l'inverse de $\mu(\ti{\sigma} )$.

\null\it Preuve: \rm Comme $\ti{\sigma} $ est unitaire, $i_{\ti{P}}^{\ti{G}}\ti{\sigma} $ l'est \'egalement, et, pour que $i_{\ti{P}}^{\ti{G}}\ti{\sigma} $ soit irr\'eductible, il faut et il suffit que $\Hom_{\ti{G}}(i_{\ti{P}}^{\ti{G}}\ti{\sigma}, i_{\ti{P}}^{\ti{G}}\ti{\sigma})$ soit de dimension $1$. Par r\'eciprocit\'e de Frobenius et le lemme g\'eom\'etrique {\bf 2.1}, ceci est certainement le cas si $w\ti{M}\ne \ti{M}$ ou $w\ti{\sigma}\ne\ti{\sigma} $.

Supposons donc maintenant $w\ti{M}=\ti{M}$ et $w\ti{\sigma} =\ti{\sigma} $. Alors, $i_{\ti{P}}^{\ti{G}}\ti{\sigma}$ est isomorphe \`a $i_{\ol{\ti{P}}}^{\ti{G}}\ti{\sigma}$. Par la r\'eciprocit\'e de Frobenius et le lemme g\'eom\'etrique {\bf 2.1}, $\Hom_{\ti{G}}(i_{\ti{P}}^{\ti{G}}\ti{\sigma}, i_{\ol{\ti{P}}}^{\ti{G}}\ti{\sigma})$ est donc de dimension $1$, si et seulement si la suite exacte $$0\rightarrow \ti{\sigma} \rightarrow (i_{\ti{P}}^{\ti{G}}\ti{\sigma} )_{\ol{\ti{P}}}\rightarrow \ti{\sigma}\rightarrow 0$$ n'est pas d\'eploy\'ee.

Remarquons que, pour $\lambda\in\Bbb R$, l'application $sp_{\lambda \alpha}:\ti{E}_{B}\rightarrow \ti{E}$, $e\otimes b\mapsto b(\chi_{\lambda\alpha })e$, d\'efinit un homomorphisme $\ti{M}$-\'equivariant $\ti{\sigma }_{B}\rightarrow\ti{\sigma }_{\lambda \alpha}$. Rappelons l'homomorphisme $q_1$  du lemme {\bf 2.2} et posons $q_{1,\lambda }=sp_{\lambda\alpha }\circ q_1$.

Si $J_{\ol{\ti{P}}\vert\ti{P}}$ n'a pas de p\^ole en $\ti{\sigma} $, alors $q_{1,0}$  se prolonge \`a tout $(i_{\ti{P}}^{\ti{G}}\ti{\sigma})_{\ol{\ti{P}}}$, et cela d\'efinit un d\'eploiement.

Inversement, si la suite exacte se d\'eploie, alors $q_{1,\lambda }$ admet un prolongement $\ti{q_{1,\lambda }}$ pour tout $\lambda\in\Bbb R$, et, alors, observant que $w\chi_{\ti{\sigma}}=\chi_{\ti{\sigma}}$ sous nos hypoth\`eses, $$\eqalign{&sp_{\lambda{\alpha } }(q_1\circ ((i_{\ti{P}}^{\ti{G}}\ti{\sigma} _{B})_{\ol{\ti{P}}}(a-\chi_{\ti{\sigma} }(a)\ ^wb_a)))\cr =&\ti{q_{1,\lambda }}\circ ((i_{\ti{P}}^{\ti{G}}\ti{\sigma} _{B})_{\ol{\ti{P}}}(a-\chi_{\ti{\sigma} }(a)\ ^wb_a)))\cr =&(\chi_{\lambda\alpha }\chi_{\ti{\sigma}})(a)-(\chi_{-\lambda\alpha }\chi_{\ti{\sigma}}) (a),\cr}$$ et on voit que cette expression s'annule en $\lambda =0$.

Donc, s'il y a un d\'eploiement, l'op\'erateur $(\chi_{\ti{\sigma} }(a)b_a- (w\chi _{\ti{\sigma} })(a) ^wb_a)J_{\ol{\ti{P}}\vert \ti{P}}$ s'annule en $\chi =1$. Mais, $\chi_{\ti{\sigma} }(a)b_a- (w\chi _{\ti{\sigma} })(a) ^wb_a$ ayant un z\'ero simple en $\chi =1$, cela signifie que $J_{\ol{\ti{P}}\vert \ti{P}}$ est r\'egulier en $\ti{\sigma} $.

Il a d\'ej\`a \'et\'e prouv\'e dans {\bf 2.2} que $J_{\ol{\ti{P}}\vert \ti{P}}$ peut au plus avoir un p\^ole simple.

Supposons maintenant $i_{\ti{P}}^{\ti{G}}\ti{\sigma} $ r\'eductible. Comme $\End_{\ti{G}}(i_{\ti{P}}^{\ti{G}}\ti{\sigma} )$ est de dimension $2$ par r\'eciprocit\'e de Frobenius, $i_{\ti{P}}^{\ti{G}}\ti{\sigma}$ doit \^etre somme directe de deux repr\'esentations non-isomorphes. L'op\'erateur d'entrelacement $\rho_w\lambda(\ti{w})J_{\ol{\ti{P}}\vert \ti{P}}(\ti{\sigma} )$ agit sur chacune d'elle par un scalaire. Il reste \`a voir que ces deux scalaires sont distincts. S'ils \'etaient \'egaux, $\rho_w\lambda(\ti{w})J_{\ol{\ti{P}}\vert \ti{P}}(\ti{\sigma} )$ serait un multiple scalaire de l'isomorphisme d'identit\'e. Autrement dit, $J_{\ol{\ti{P}}\vert \ti{P}}(\ti{\sigma} )$ serait un multiple scalaire de $\rho_w^{-1}\lambda(\ti{w}^{-1})$. Or, par r\'eciproci-t\'e de Frobenius, $\rho_w^{-1}\lambda(\ti{w}^{-1})$ correspond \`a l'homomorphisme $(i_{\ti{P}}^{\ti{G}}\ti{\sigma} )_{\ol{\ti{P}}}\rightarrow\ti{\sigma }$ induite par $v\mapsto \rho_w^{-1}(v(\ti{w}^{-1}))$. Cet homomorphisme est nul sur le sous-espace des \'el\'ements de $i_{\ti{P}}^{\ti{G}}\ti{E}$ \`a support dans $\ti{P}\ol{\ti{P}}$. De l'autre c\^ot\'e, par le lemme {\bf 2.2}, $J_{\ol{\ti{P}}\vert \ti{P}}(\ti{\sigma} )$ vient par r\'eciprocit\'e de Frobenius d'un homomorphisme qui n'est pas identiquement nul sur ce sous-espace.

Composons l'automorphisme $\rho_w\lambda(\ti{w})J_{\ol{\ti{P}}\vert \ti{P}}(\ti{\sigma} )$ avec lui-m\^eme, on trouve $$\rho_w\lambda(\ti{w})\rho_w\lambda(\ti{w})J_{\ti{P}\vert \ol{\ti{P}}}(\ti{\sigma} )J_{\ol{\ti{P}}\vert \ti{P}}(\ti{\sigma} )=\rho_w^2\lambda(\ti{w}^2) \gamma(\ti{G}/\ti{M}) \mu(\ti{\sigma })^{-1}.$$ On observe que $\rho_w^2\lambda(\ti{w}^2)$ est un automorphisme de $i_{\ti{P}}^{\ti{G}}\ti{\sigma} $ qui vient par fonctorialit\'e d'un automorphisme de $\ti{\sigma }$. Il est donc scalaire, et on peut normaliser $\rho _w$ de telle sorte que ce scalaire soit $\gamma(\ti{G}/\ti{M})^{-1}$. Dans ce cas, $\rho_w\lambda(\ti{w})J_{\ol{\ti{P}}\vert \ti{P}}(\ti{\sigma} )$ agit donc sur les deux sous-espaces par des scalaires de signe oppos\'e de carr\'e \'egal \`a $\mu(\ti{\sigma })^{-1}$. Ceci termine la preuve de la proposition. \hfill{\fin 2}

\null{\bf 3.2 Lemme:} \it Soit $\lambda >0$ tel que $i_{\ti{P}}^{\ti{G}}\ti{\sigma}_{\lambda\alpha }$ soit r\'eductible.

Alors, $w\ti{M}=\ti{M}$ et $w\ti{\sigma} =\ti{\sigma} $.

La repr\'esentation induite $i_{\ti{P}}^{\ti{G}}\ti{\sigma}_{\lambda\alpha }$ est de longueur $2$.

Elle a une unique sous-repr\'esentation irr\'eductible qui est de carr\'e int\'egrable et un unique quotient irr\'eductible. Les deux ne sont pas isomorphes.

La repr\'esentation $i_{\ti{P}}^{\ti{G}}\ti{\sigma}_{-\lambda\alpha }$ est \'egalement r\'eductible (et vice-versa sa r\'eductibilit\'e implique celle de $i_{\ti{P}}^{\ti{G}}\ti{\sigma}_{\lambda\alpha }$). Son quotient irr\'eductible est la sous-repr\'esentation irr\'e-ductible de $i_{\ti{P}}^{\ti{G}}\ti{\sigma}_{\lambda\alpha }$, et le quotient irr\'eductible de $i_{\ti{P}}^{\ti{G}}\ti{\sigma}_{\lambda\alpha }$  est la sous-repr\'esentation  irr\'eductible de $i_{\ti{P}}^{\ti{G}}\ti{\sigma}_{-\lambda\alpha }$.

En particulier, $i_{\ti{P}}^{\ti{G}}\ti{\sigma} _{\lambda\alpha }$ n'est pas semi-simple et l'op\'erateur $J_{\ol{\ti{P}}\vert \ti{P}}(\ti{\sigma} _{\lambda\alpha })$ est bien d\'efini, non nul, mais pas bijectif.

\null Preuve: \rm Si $w\ti{M}\ne \ti{M}$, alors l'irr\'eductibilit\'e de $i_{\ti{P}}^{\ti{G}}\ti{\sigma}_{\lambda\alpha }$ dans le cas $\widetilde{G}$ r\'eductif est prouv\'e dans \cite{Cs, 7.1}. On ne va pas recopier cette preuve ici dont les id\'ees sont semblables \`a ce qui suivra dans le cas o\`u cette induite est r\'eductible, mais cette preuve n'utilise que le lemme g\'eom\'etrique {\bf 2.1}, les propri\'et\'es du foncteur de Jacquet, et le crit\`ere de Casselman pour les repr\'esentations de carr\'e int\'egrable qui se g\'en\'eralisent tous au cas d'une extension centrale.

Si $i_{\ti{P}}^{\ti{G}}\ti{\sigma}_{\lambda \alpha}$ est r\'eductible, par le lemme g\'eom\'etrique {\bf 2.1}, il est n\'ecessairement de longueur $2$, et $(i_{\ti{P}}^{\ti{G}}\ti{\sigma}_{\lambda \alpha})_P$ a deux sous-quotients irr\'eductibles cuspidaux non-isomorphes $\ti{\sigma}_{\lambda\alpha}$ et $(w\ti{\sigma})_{-\lambda\alpha }$. La r\'eciprocit\'e de Frobenius montre alors que $i_{\ti{P}}^{\ti{G}}\ti{\sigma}_{\lambda\alpha }$ ne peut pas \^etre semisimple et par le crit\`ere de Casselman l'unique sous-repr\'esentation irr\'eductible $\ti{\pi _1}$ est de carr\'e int\'egrable. Notons $\ti{\pi _2}$ l'unique quotient irr\'eductible. Alors, $(\ti{\pi _1})_{\ti{P}}=\ti{\sigma}_{\lambda\alpha }$, $(\ti{\pi _2})_{\ti{P}}=w(\ti{\sigma}_{\lambda\alpha })$ et on a une suite exacte $0\rightarrow \ti{\pi _1}\rightarrow i_{\ti{P}}^{\ti{G}}\ti{\sigma}_{\lambda\alpha }\rightarrow\ti{\pi _2}\rightarrow 0$. Comme $\ti{\pi _1}$ est de carr\'e int\'egrable, la repr\'esentation conjugu\'ee $\ol{\ti{\pi _1}}$ est isomorphe \`a la repr\'esentation contragr\'ediente $\ti{\pi _1}^{\vee }$ de $\ti{\pi _1}$.

On a $(\ol{\ti{\pi _1}})_{\ti{P}}=\ol{\ti{\sigma} }_{\lambda\alpha }$ et $(\ti{\pi _1}^{\vee })_{\ti{P}}\simeq(\ti {\pi _1})_{\ol{\ti{P}}}^{\vee }=(w(\ti{\sigma}_{\lambda\alpha }))^{\vee }\simeq(w\ol{\ti{\sigma} })_{\lambda\alpha }$ (voir aussi \cite{KS, preuve de 3.10}). Comme les deux sont isomorphes, on en d\'eduit $\ol{\ti{\sigma}}\simeq w\ol{\ti{\sigma} }$ et donc $\ti{\sigma}\simeq w\ti{\sigma}$.

Par cons\'equence, $(i_{\ti{P}}^{\ti{G}}\ti{\sigma}_{\lambda\alpha })_{\ti{P}}\simeq\ti{\sigma}_{\lambda\alpha }\oplus \ti{\sigma}_{-\lambda\alpha }$ et $(\ti{\pi _2})_{\ti{P}}=\ti{\sigma}_{-\lambda\alpha }$. Il r\'esulte alors du lemme g\'eom\'etrique {\bf 2.1} et de la r\'eciprocit\'e de Frobenius que $\ti{\pi _2}$ est la sous-repr\'esentation  irr\'eductible de $i_{\ti{P}}^{\ti{G}}\ti{\sigma}_{-\lambda\alpha }$ et que $\ti{\pi _1}$ en est le quotient irr\'eductible. La repr\'esentation $i_{\ti{P}}^{\ti{G}}\ti{\sigma}_{-\lambda\alpha }$ n'est donc pas isomorphe \`a $i_{\ti{P}}^{\ti{G}}\ti{\sigma}_{\lambda\alpha }$, et ainsi $J_{\ol{\ti{P}}\vert \ti{P}}(\ti{\sigma} _{\lambda\alpha })$ est r\'egulier et non nul suite au corollaire {\bf 2.3}, mais il ne peut pas \^etre bijectif.
\hfill{\fin 2}

\null{\bf 3.3} Pour la preuve du lemme suivant, il sera commode de remplacer l'anneau des fonctions $B$ du tore complexe $\X(\ti{M})$ par celui correspondant \`a sa sous-vari\'et\'e $\{\ti{\chi }_{\lambda\alpha }\vert\lambda\in\Bbb C\}$: observons d'abord que $\ti{\chi}_{\alpha }$ d\'efinit un homomorphisme $M/M^1\rightarrow q^{\Bbb Z}$. Choisissons $m_0\in\ti{M}$ tel que $\ti{\chi}_{\alpha }(m_0)$ engendre l'image de cet homomorphisme. D\'efinissons pour tout $m\in \ti{M}$ un \'el\'ement $m'\in m_0^{\Bbb Z}$ par $\ti{\chi}_{\alpha }(m)=\ti{\chi}_{\alpha }(m')$. Il est d\'etermin\'e de fa\c con unique. L'application $m\mapsto m'$ ainsi d\'efinie est un homomorphisme de groupes. Posons $Z=b_{m_0}$, notons $B^1$ l'anneau des polyn\^omes de Laurent $\Bbb C[Z,Z^{-1}]$, $\ti{E}_{B^1}=\ti{E}\otimes B^1$ et $\ti{\sigma }_{B^1}$ la repr\'esentation de $\ti{M}$ dans $\ti{E}_{B^1}$ donn\'ee par $(\ti{\sigma }_{B^1}(m))(e\otimes b)=\ti{\sigma}(m)e\otimes bb_{m'}$. Pour $\lambda\in\Bbb C$, \'ecrivons $z_{\lambda }=\ti{\chi}_{\lambda\alpha }(m_0)$ et notons ${\frak m}_{\lambda }$ l'id\'eal principal de $B^1$ engendr\'e par $Z-z_{\lambda }$. C'est un id\'eal maximal et l'application canonique $sp_{\lambda }:\ti{E}_{B^1}\rightarrow \ti{E}\otimes B^1/{\frak m}_{\lambda }$ d\'efinit un homomorphisme entre les repr\'esentations $\ti{\sigma }_{B^1}$ et $\ti{\sigma }_{\lambda\alpha}$.

Remarquons que, si $G$ est un groupe r\'eductif semi-simple, $B^1=B$ et $Z=Y$. En g\'en\'eral, $Y$ est un mon\^ome en $Z$, en particulier la fonction $\mu $ et l'op\'erateur d'entrelacement $J_{\ol{P}\vert P}$ sont rationnelles en $B^1$.

\null {\bf Lemme:} \it Les p\^oles de la fonction $\mu $ de Harish-Chandra sont simples.

\null Preuve: \rm On a vu en {\bf 2.3} que $J_{\ol{\ti{P}}\vert \ti{P}}$ est r\'egulier en $\ti{\sigma}_{\lambda\alpha}$, si $\lambda\ne 0$, et on sait que $\mu $ est r\'eguli\`ere en $\ti{\sigma}$ (cf. \cite{W, V.2.1}, \cite{Li1, 2.4.3}). Donc, si $\mu $ a un p\^ole en $\ti{\sigma}_{\lambda\alpha}$, alors $\lambda \ne 0$, et au moins un des deux op\'erateurs, $J_{\ol{\ti{P}}\vert \ti{P}}(\ti{\sigma}_{\lambda\alpha})$ ou $J_{\ti{P}\vert \ol{\ti{P}}}(\ti{\sigma}_{\lambda\alpha})$ n'est pas bijectif. Mais, comme aucun des deux ne peut \^etre nul (cf. {\bf 2.3}), les deux ne sont pas bijectifs. Donc, $i_{\ti{P}}^{\ti{G}}\ti{\sigma}_{\lambda\alpha }$ est r\'eductible, ce qui implique $w\ti{\sigma}\simeq \ti{\sigma}$ par {\bf 3.2}.
Par ailleurs, par sym\'etrie \cite{W, V.2.1}, il suffira de montrer que $\mu $ a un p\^ole au plus simple en $\ti{\sigma}_{\lambda\alpha}$, si $\lambda >0$.

Soit donc $\lambda >0$ tel que $\mu $ ait un p\^ole en $\ti{\sigma}_{\lambda\alpha}$. Pour simplifier les notations dans la suite, posons $\ti{\tau}=\ti{\sigma}_{\lambda\alpha}$. Il suffira de montrer que l'op\'erateur compos\'e $J_{\ol{\ti{P}}\vert \ti{P}}J_{\ti{P}\vert \ol{\ti{P}}}$ a un z\'ero simple en $\ti{\tau}$.

Par {\bf 3.2}, $i_{\ti{P}}^{\ti{G}}\ti{\tau}$ est r\'eductible, n\'ecessairement de longueur $2$ et non semi-simple.  On a une suite exacte $0\rightarrow\ti{\pi_1}\rightarrow i_{\ti{P}}^{\ti{G}}\ti{\tau }\rightarrow\ti{\pi_2}\rightarrow 0$, o\`u $\ti{\pi _1}$ et $\ti{\pi _2}$ sont des repr\'esentations irr\'eductibles et que $\ti{\pi _1}$ est de carr\'e int\'egrable.

Posons $\ti{\Pi}=i_{\ti{P}}^{\ti{G}}\ti{\sigma}_{B^1}\otimes_{B^1} B^1/{\frak m}_{\lambda }^2$. Notons $\ti{V}$ l'espace de $\ti{\Pi}$. Remarquons que la repr\'esentation $\ti{\Pi}$ est admissible, puisque de longueur finie.

Comme l'id\'eal ${\frak m}_{\lambda }$ est principal, on a une suite exacte $$0\rightarrow i_{\ti{P}}^{\ti{G}}\ti{\sigma }_{B^1}\otimes_{B^1} {\frak m}_{\lambda }/{\frak m}_{\lambda }^2 \rightarrow \ti{\Pi }\rightarrow i_{\ti{P}}^{\ti{G}}\ti{\sigma }_{B^1}\otimes_{B^1} B^1/{\frak m}_{\lambda }\rightarrow 0,$$ o\`u le quotient et la sous-repr\'esentation sont tous les deux isomorphes \`a $i_{\ti{P}}^{\ti{G}}\ti{\tau }$. On en d\'eduit une filtration de Jordan-H\"older de $\ti{\Pi }$ de longueur $4$ donn\'ee par
$$0\subsetneq \ti{V_1}\subsetneq\ti{V_2}\subsetneq\ti{V_3}\subsetneq\ti{V_4}=\ti{V},$$ o\`u $\ti{V_1}\simeq\ti{\pi _1}$, $\ti{V_2}/\ti{V_1}\simeq\ti{\pi _2}$, $\ti{V_3}/\ti{V_2}\simeq\ti{\pi _1}$ et $\ti{V_4}/\ti{V_3}\simeq\ti{\pi_2}$.

On montrera tout-\`a-l'heure que toute d\'ecomposition de Jordan-H\"older  de $\ti{V}$ est de cette forme. Ceci implique l'\'enonc\'e: en effet, les $\ti{V_i}$ seront alors les seuls sous-modules invariants de $\ti{G}$, et on aura de m\^eme une unique filtration pour $\ti{\Pi'}=i_{\ol{\ti{P}}}^{\ti{G}}\sigma_{B^1}\otimes_{B^1} B^1/{\frak m}_{\lambda }^2$ que l'on notera $0\subsetneq\ti{V'_1}\subsetneq \ti{V'_2}\subsetneq\ti{V'_3}\subsetneq\ti{V'_4}=\ti{V'},$ $\ti{V'}$ d\'esignant l'espace de $\ti{\Pi '}$, avec les sous-quotients $\ti{V_1'}\simeq\ti{\pi_2}$, $\ti{V_2'}/\ti{V_1'}\simeq\ti{\pi _1}$, $\ti{V_3'}/\ti{V_2'}\simeq\ti{\pi _2}$ et $\ti{V_4'}/\ti{V_3'}\simeq\ti{\pi_1}$ par {\bf 3.2}.

L'op\'erateur d'entrelacement $J_{\ol{\ti{P}}\vert \ti{P}}(\ti{\sigma }_{B^1})$ tensoris\'e avec l'identit\'e sur $B^1/{\frak m}_{\lambda}^2$ envoie $\ti{V}$ dans $\ti{V'}$ et, comme $J_{\ol{\ti{P}}\vert \ti{P}}(\ti{\tau })$ est non nul (cf. {\bf 2.3}), il doit envoyer $\ti{V_1}$ sur $0$, $\ti{V_2}$ sur $\ti{V_1'}$, $\ti{V_3}$ sur $\ti{V_2'}$ et $\ti{V_4}$ sur $\ti{V_3'}$. De m\^eme, l'op\'erateur d'entrelacement $J_{\ti{P}\vert \ol{\ti{P}}}(\ti{\sigma }_{B^1})$ tensoris\'e avec l'identit\'e sur $B^1/{\frak m}_{\lambda}^2$ envoie $\ti{V'}$ dans $\ti{V}$ et plus pr\'ecis\'ement  $\ti{V'_1}$ sur $0$, $\ti{V'_2}$ sur $\ti{V_1}$, $\ti{V'_3}$ sur $\ti{V_2}$ et $\ti{V_4'}$ sur $\ti{V_3}$. On en d\'eduit que l'op\'erateur compos\'e $J_{\ti{P}\vert \ol{\ti{P}}}(\ti{\sigma }_{B^1})J_{\ol{\ti{P}}\vert \ti{P}}(\ti{\sigma }_{B^1})$ tensoris\'e avec l'identit\'e sur $B^1/{\frak m}_{\lambda}^2$ envoie $\ti{V_3}$ sur $\ti{V_1}$ et $\ti{V_4}$ sur $\ti{V_2}$. En particulier, il est non nul. Ceci montre que le z\'ero en $\ti{\tau }=\ti{\sigma}_{\lambda\alpha }$ est simple: en effet, sinon les coefficients de $J_{\ti{P}\vert \ol{\ti{P}}}(\ti{\sigma }_{B^1})J_{\ol{\ti{P}}\vert \ti{P}}(\ti{\sigma }_{B^1})$  serait tous divisibles par $(X-X(\chi_{\lambda\alpha}))^2$, donc par $(Z-z_{\lambda})^2$, puisque $X$ est une puissance de $Z$. Ce seraient donc des \'el\'ements de ${\frak m}_{\lambda}^2$ et, compte tenu du produit tensoriel avec $B^1/{\frak m}_{\lambda}^2$, l'op\'erateur tensoris\'e serait effectivement $0$.

Reste \`a montrer que toute d\'ecomposition de $\ti{V}$ a la forme indiqu\'ee. Supposons par absurde qu'il existe un sous-espace invariant $\ti{U}$ de $\ti{V}$, $0\ne\ti{U}\ne\ti{V}$, diff\'erent de $\ti{V_1}$, $\ti{V_2}$ et $\ti{V_3}$. Alors, $\ti{U}\cap \ti{V_2}\in\{0,\ti{V_1},\ti{V_2}\}$.

Si $\ti{U}\cap\ti{V_2}=\ti{V_2}$, alors ou bien $\ti{U}=\ti{V_2}$ ou bien $\ti{U}/\ti{V_2}$ est un sous-espace non trivial de $\ti{V}/\ti{V_2}$, donc n\'ecessairement \'egal \`a $\ti{V_3}/\ti{V_2}$, i.e. $\ti{U}=\ti{V_3}$. De m\^eme, si $\ti{U}\subsetneq\ti{V_2}$, on a n\'ecessairement $\ti{U}=\ti{V_1}$.

Donc, $\ti{V_2}\not\subseteq\ti{U}$, et $\ti{V_2}+\ti{U}\in\{\ti{V_3},\ti{V_4}\}$. Il existe donc $\ti{U'}\subseteq\ti{U}$ tel que $\ti{V_2}+\ti{U'}=\ti{V_3}$. Regardons son intersection avec $\ti{V_2}$. Il appartient \`a $\{0,\ti{V_1}, \ti{V_2}\}$. Si c'\'etait $\ti{V_2}$, alors $\ti{V_2}\subseteq\ti{U'}\subseteq \ti{U}$, ce qui avait d\'ej\`a \'et\'e exclu. Si $\ti{U'}\cap\ti{V_2}=0$, alors $\ti{V_3}=\ti{V_2}\oplus\ti{U'}$. Posons alors $\ti{U''}=\ti{U'}\oplus\ti{V_1}$. C'est un sous-espace semisimple isomorphe \`a $\ti{\pi_1}\oplus\ti{\pi_1}$. Finalement, si $\ti{U'}\cap\ti{V_2}=\ti{V_1}$, alors $\ti{U'}/\ti{V_1}\simeq\ti{\pi _1}$. On voit donc que $\ti{V}$ contient dans tous les cas une sous-repr\'esentation $\ti{U''}$ dont la semisimplifi\'ee est $\ti{\pi_1}\oplus\ti{\pi _1}$. On va montrer que $\ti{\Pi }$ ne peut pas poss\'eder une telle sous-repr\'esentation. Ceci donnera la contradiction souhait\'ee.

Observons d'abord que $\ti{\Pi }$ poss\`ede un caract\`ere central unitaire \'egal \`a la restriction de celui de $\ti{\sigma }$ au centre de $\ti{G}$: en effet, le centre de $\ti{G}$ agit trivialement sur $B^1$. Il agit donc par la repr\'esentation induite $i_{\ti{P}}^{\ti{G}}\ti{\sigma }$ dont le caract\`ere central a la forme indiqu\'ee. Comme $\ti{\pi _1}$ est de carr\'e int\'egrable, $\ti{\Pi }_{\vert\ti{U''}}$ v\'erifie alors le crit\`ere de Casselman (voir par exemple \cite{W, III.1}). On conclut que $\ti{\Pi }_{\vert\ti{U''}}$ est lui-m\^eme de carr\'e int\'egrable et en particulier unitaire, i.e. $\ti{\Pi }_{\vert\ti{U''}}=\ti{\pi_1}\oplus\ti{\pi _1}$.

Il reste \`a voir que $\ti{\Pi }$ ne peut pas poss\'eder une telle sous-repr\'esentation. Par r\'eciprocit\'e de Frobenius, $r_{\ti{P}}^{\ti{G}}\ti{\pi _1}=\ti{\tau}$ (voir preuve du lemme {\bf 3.2}), alors que, par le lemme g\'eom\'etrique \cite{BZ, 5.2}, $r_{\ti{P}}^{\ti{G}}\ti{\Pi}=(\ti{\sigma}_{B^1}\otimes_{B^1} B^1/{\frak m}_{\lambda }^2)\oplus w(\ti{\sigma}_{B^1}\otimes_{B^1} B^1/{\frak m}_{\lambda }^2)$. Ici, on a d\'ecompos\'e $r_{\ti{P}}^{\ti{G}}\ti{\Pi}$ en somme directe selon ses exposants pour $\ti{A^\dagger}$ (cf. \cite{W, I.3}), utilisant le fait que les sous-quotients du terme \`a gauche sont tous isomorphes \`a $\ti{\sigma }_{\lambda\alpha }$ alors que ceux du terme \`a droite sont tous isomorphes \`a $\ti{\sigma }_{-\lambda\alpha }$. On en d\'eduirait que  $\ti{\tau}\oplus\ti{\tau}\simeq \ti{\sigma}_{B^1}\otimes_{B^1} B^1/{\frak m}_{\lambda }^2$ par comparaison des exposants.

Or, ceci est absurde: le groupe $\ti{A^\dagger}$ agit sur l'espace de $\ti{\tau}\oplus\ti{\tau }$ par le caract\`ere central de $\ti{\tau }$ qui est celui de $\ti{\sigma }$, alors que l'action de $\ti{A^\dagger}$ sur $\ti{\sigma}_{B^1}\otimes_{B^1} B^1/{\frak m}_{\lambda }^2$ n'est pas semi-simple: elle se fait par $a\mapsto\chi_{\sigma }(a)(1+b_{a'})$, et l'action de $b_{a'}$ n'est pas semi-simple si $a$ n'est pas dans le centre de $\ti{G}$: en effet, $b_{a'}$ est alors de la forme $Z^k$ pour un certain entier $k\ne 0$ et $Z^k$ n'est jamais congru \`a un nombre complexe modulo $(Z-z_{\lambda })^2$.
\hfill{\fin 2}

\null{\bf 3.4 Lemme:} \it (i) Si $w\ti{\sigma }\simeq\ti{\sigma }$, l'application $\lambda\mapsto\mu(\ti{\sigma}_{\lambda\alpha })$ est \`a valeurs r\'eelles pour tout $\lambda $ r\'eel en lequel cette application est d\'efinie.

(ii) Si $\mu (\ti{\sigma} )=0$, alors $w\ti{M}=\ti{M}$, $w\ti{\sigma }\simeq\ti{\sigma }$ et l'application $\lambda\mapsto\mu(\ti{\sigma}_{\lambda\alpha })$ a un unique p\^ole $\lambda _0>0$. Il v\'erifie par ailleurs $\lambda_0\alpha\leq\delta _{\ti{P}}^{1/2}$. L'unique p\^ole sur l'axe n\'egatif est alors $-\lambda _0$.

(iii) Si $\mu(\ti{\sigma} )\ne 0$, alors $\mu(\ti{\sigma} )>0$ et l'application $\lambda\mapsto\mu(\ti{\sigma}_{\lambda\alpha })$ est r\'eguli\`ere et $>0$ pour tout $\lambda\in\Bbb R$.

\null Preuve: \rm Observons d'abord que, pour $\lambda\in\Bbb R$, la contragr\'ediente de $\ti{\sigma} $ \'etant isomorphe \`a sa conjugu\'ee et la fonction $\mu $ invariante par $w$ et passage \`a la contragr\'ediente \cite{W, V.2.1} et \cite{Li1, 2.4.3},
$$\ol{\mu(\ti{\sigma}_{\lambda\alpha})}=\mu(\ol{\ti{\sigma} }_{\lambda\alpha})=\mu((\ti{\sigma}^{\vee })_{\lambda\alpha })=\mu(\ti{\sigma}_{-\lambda_{\alpha }})=\mu(w(\ti{\sigma} _{\lambda\alpha }))=\mu(\ti{\sigma}_{\lambda\alpha}),$$ ce qui prouve $\mu(\ti{\sigma}_{\lambda\alpha})\in\Bbb R$.

Supposons maintenant $\mu (\ti{\sigma} )=0$. Par {\bf 2.3}, on a alors $w\ti{M}=\ti{M}$ et $w\ti{\sigma }\simeq\ti{\sigma }$. S'il n'y avait aucun p\^ole sur la demi-droite $\lambda>0$, comme il n'y a pas non plus de z\'ero pour $\lambda\ne 0$, l'op\'erateur d'entrelacement normalis\'e $R_{\ti{P}}(\ti{\sigma} _{\lambda\alpha },\ti{w})$ de {\bf 2.6} serait bijectif pour tout $\lambda\geq 0$. Comme il est positif d\'efini pour $\lambda =0$, il resterait positif d\'efini pour tout $\lambda\geq 0$, ce qui impliquerait que la repr\'esentation $i_{\ti{P}}^{\ti{G}}\ti{\sigma}_{\lambda\alpha }$ serait irr\'eductible unitaire pour tout $\lambda\geq 0$.

Montrons que ceci est impossible et qu'en fait il doit y avoir un p\^ole. En effet, si ces repr\'esentations \'etaient toutes unitaires, par l'in\'egalit\'e de Cauchy-Schwarz, les coefficients matriciels de $i_{\ti{P}}^{\ti{G}}\ti{\sigma}_{\lambda\alpha}$ seraient born\'es pour tout $\lambda>0$. Or, par la formule d'approximation de Casselman (cf. {\bf 2.5}), ceci est impossible lorsque $\lambda\alpha >\delta_{\ti{P}}^{1/2}$, d'o\`u l'existence d'un p\^ole $\lambda_0$ qui v\'erifie par ailleurs $\lambda_0\alpha\leq\delta _{\ti{P}}^{1/2}$.

Quant \`a l'unicit\'e, rappelons que, pour $\lambda>0$, l'unique sous-repr\'esentation de $i_{\ti{P}}^{\ti{G}}\ti{\sigma}_{\lambda\alpha}$ est par {\bf 3.2} de carr\'e int\'egrable, si cette repr\'esentation est r\'eductible, ce qui est le cas par {\bf 2.4} (i) si $\lambda $ est un p\^ole de $\mu(\ti{\sigma}_{\cdot\alpha})$. On sait par {\bf 3.3} que ce p\^ole est simple et par {\bf 2.4} (ii) que son r\'esidu est $>0$. On observe que $\mu $ a un z\'ero d'ordre $2$ en $\ti{\sigma} $, puisque les deux op\'erateurs d'entrelacement qui d\'efinissent $\mu $ ont alors chacun un p\^ole simple en $\ti{\sigma }$ par {\bf 3.1}. De plus, les valeurs de $\mu(\ti{\sigma}_{\lambda\alpha})$ sont r\'elles pour $\lambda $ r\'eel par la partie (i). Comme le premier r\'esidu est $>0$, on en d\'eduit que $\mu(\ti{\sigma}_{\lambda\alpha})<0$ sur l'axe r\'eel au voisinage de $0$ \`a l'exclusion de $\lambda=0$. Par contre, apr\`es ce r\'esidu, vu que le p\^ole est simple,  on a $\mu(\ti{\sigma}_{\lambda\alpha})>0$ jusqu'au prochain p\^ole ou prochain z\'ero. Or, il n'y a plus de z\'ero et s'il y avait un p\^ole qui serait n\'ecessairement simple, on aurait un changement de signe du positif vers le n\'egatif. Le r\'esidu serait donc $<0$. Or, ceci n'est pas possible par {\bf 2.4} (ii). Donc, il n'y a qu'un seul p\^ole.

Reste \`a consid\'erer le cas $\mu(\ti{\sigma} )\ne 0$. Alors, $\mu(\ti{\sigma} )>0$ par \cite{W, V.2.1} et \cite{Li1, 2.4.3}. Comme $\lambda\mapsto\mu(\ti{\sigma}_{\lambda\alpha })$ ne s'annule pas pour $\lambda>0$, les valeurs de cette application restent $>0$ jusqu'au premier p\^ole \'eventuel. Or, celui-ci \'etant n\'ecessairement simple, on voit comme ci-dessus que son r\'esidu serait $<0$, ce qui contredit {\bf 2.4}. Par sym\'etrie, on conclut que $\mu(\ti{\sigma}_{\lambda\alpha })>0$ pour tout $\lambda\in\Bbb R$.\hfill{\fin 2}

\null{\bf 3.5} Le r\'esultat suivant et sa preuve sont bien connus dans le cas $p$-adique. On verra qu'ils restent valable dans le cas d'une extension centrale.

\null{\bf Corollaire:} \it La repr\'esentation $i_{\ti{P}}^{\ti{G}}\ti{\sigma}_{\lambda\alpha }$ est unitaire pour $\lambda =0$. Si $\mu (\ti{\sigma} )\ne 0$, elle est irr\'eductible et non unitaire pour $\lambda\ne 0$. Si $\mu(\ti{\sigma} )=0$, elle est unitaire jusqu'au premier p\^ole de $\mu $.

\null Preuve: \rm Si $w\ti{\sigma} \not\simeq\ti{\sigma}$, il r\'esulte de la r\'eciprocit\'e de Frobenius et du lemme g\'eom\'etrique \cite{BZ, 5.2} que $i_{\ti{P}}^{\ti{G}}\ti{\sigma}_{\lambda\alpha }$ n'est pas hermitienne, donc pas non plus unitaire pour $\lambda\ne 0$. Supposons maintenant $\ti{\sigma}\simeq w\ti{\sigma} $. Si $\mu(\ti{\sigma} )\ne 0$, alors l'op\'erateur d'entrelacement normalis\'e $R_{\ti{P}}(\ti{\sigma},\ti{w})$ de {\bf 2.6} agit par {\bf 3.1} par des scalaires de signe diff\'erent et cela sera ainsi jusqu'au premier p\^ole de $\mu (\sigma_{\cdot\alpha })$. Or, il n'y a pas de tel p\^ole par {\bf 3.4} (iii). Si $\mu(\ti{\sigma} )=0$, alors il y a un p\^ole $\lambda_0$ par {\bf 3.4} (ii), et $R_{\ti{P}}(\ti{\sigma},\ti{w})$ s'annule sur l'unique sous-repr\'esentation irr\'eductible de $i_{\ti{P}}^{\ti{G}}\ti{\sigma}_{\lambda_0\alpha }$. Ce z\'ero est  n\'ecessairement d'ordre $1$ par {\bf 3.3}. Or, la signature de $R_{\ol{\ti{P}}\vert \ti{P}}(\ti{\sigma} )$ sur ce sous-espace devient alors $<0$ pour $\lambda>\lambda_0$ et reste $>0$ sur un compl\'ement de ce sous-espace. Comme il n'y a plus d'autres p\^oles ou z\'eros de $\mu (\sigma_{\cdot\alpha })$, la repr\'esentation ne peut plus redevenir unitaire.

\null{\bf 4.} Dans cette section, on fera la preuve du th\'eor\`eme {\bf 1.1} qui sera un assemblage des lemmes prouv\'es en section {\bf 3.} ainsi que des r\'esultats rappel\'es en section {\bf 2.}:

\null \it Preuve du th\'eor\`eme:\rm

Montrons d'abord la propri\'et\'e (i). Il r\'esulte du lemme {\bf 3.2} (resp. lemme {\bf 3.1}), pour $\lambda\ne 0$ (resp. $\lambda=0$), que l'induite est irr\'eductible si $w\ti{M}\ne \ti{M}$ ou si $w\ti{\sigma}\not\simeq\ti{\sigma} $. Le corollaire {\bf 2.3} dit que l'op\'erateur d'entrelacement est r\'egulier sous les hypoth\`eses de (i), donc la fonction $\mu $ est non nulle. Ceci montre \'egalement que, si $\mu $ avait un p\^ole, alors au moins un des deux op\'erateurs d'entrelacement ne serait pas bijectif, ce qui contredit l'irr\'eductibilit\'e de l'induite.

Quant \`a (ii), montrons d'abord b) et c): l'irr\'eductibilit\'e de $i_{\ti{P}}^{\ti{G}}\ti{\sigma} $ dans b) est \'equivalent par le lemme {\bf 3.1}, \`a ce que l'op\'erateur d'entrelacement ait un p\^ole en $\ti{\sigma}$. Or, l'autre op\'erateur d'entrelacement devant alors \'egalement avoir un p\^ole, ceci \'equivaut \`a $\mu(\ti{\sigma} )=0$. La deuxi\`eme partie de b) correspond \`a la deuxi\`eme partie du lemme {\bf 3.1}. Quant \`a c), la simplicit\'e du p\^ole est l'objet du lemme {\bf 3.3} et le fait que la fonction inverse de $\mu $ est d\'efinie en $\ti{\sigma }_{\lambda\alpha }$ r\'esulte du lemme {\bf 2.2} qui dit que les deux op\'erateurs d'entrelacement qui le composent sont dans ce cas r\'egulier en $\ti{\sigma }_{\lambda\alpha }$. On a d\'ej\`a v\'erifi\'e dans la proposition {\bf 2.4}(i) que l'induite est r\'eductible, si $\mu $ a un p\^ole en $\lambda >0$. R\'eciproquement, s'il y a r\'eductibilit\'e, l'induite ne peut pas \^etre semisimple par {\bf 3.2}, et il n'est donc pas isomorphe \`a l'induite \`a partir de $\ol{\ti{P}}$, les repr\'esentations et quotients irr\'eductibles n'\'etant pas les m\^emes. Donc, les deux op\'erateurs d'entrelacement qui sont r\'eguliers par le corollaire {\bf 2.3} ne sont pas bijectifs. Leur compos\'ee est donc identiquement nulle, i.e. la fonction $\mu $ a un p\^ole en ce $\lambda>0$.

Passons maintenant au a): si l'induite est irr\'eductible pour $\lambda =0$, alors on sait par b) que $\mu(\ti{\sigma} )=0$. Il existe alors un unique p\^ole $\lambda _0>0$ de la fonction $\mu$ par le lemme {\bf 3.4} (ii). La r\'eductibilit\'e de l'induite r\'esulte de c) comme l'unicit\'e du point de r\'eductibilit\'e $>0$, compte tenu de l'unicit\'e du p\^ole que l'on vient de remarquer. Reste \`a montrer l'unicit\'e du point de r\'eductibilit\'e en cas de r\'eductibilit\'e en $\lambda =0$. Or, par le lemme {\bf 3.4} (iii), la fonction $\mu $ n'a alors pas de p\^ole sur l'axe r\'eelle et donc par c) il n'y a pas non plus de point de r\'eductibilit\'e correspondant \`a un $\lambda\ne 0$.

(iii) Les p\^oles de $J_{\ol{\ti{P}}\vert \ti{P}}$ sont d'ordre $1$ par le lemme {\bf 3.1} ou {\bf 2.3}. Le lemme {\bf 2.2} montre que, si $J_{\ol{\ti{P}}\vert \ti{P}}$ a un p\^ole en $\ti{\sigma}_{\lambda\alpha}$, $\lambda\in\Bbb R$, alors $\lambda=0$ et $i_{\ti{P}}^{\ti{G}}\ti{\sigma }$ est irr\'eductible par {\bf 3.1}, ce qui implique $\mu(\ti{\sigma} )=0$ par le ii) b) du th\'eor\`eme. R\'eciproquement, si $\ti{\sigma}_{\lambda\alpha}$, $\lambda\in\Bbb R$, est un z\'ero de $\mu $, alors, par le (i) du th\'eor\`eme, $w\ti{M}=\ti{M}$ et $w\ti{\sigma }=\ti{\sigma }$ et, par le (ii) c) du th\'eor\`eme, $\lambda =0$. On d\'eduit alors du (ii) b) du th\'eor\`eme que  $i_{\ti{P}}^{\ti{G}}\ti{\sigma }$ est irr\'eductible et ensuite du lemme {\bf 3.1} que l'op\'erateur d'entrelacement $J_{\ol{\ti{P}}\vert \ti{P}}$ a un p\^ole en $\ti{\sigma}$.

(iv) Observons d'abord que, si $w\ti{M}\ne \ti{M}$ ou $w\ti{\sigma }\ne\ti{\sigma }$, alors par la partie a) du th\'eor\`eme $i_{\ti{P}}^{\ti{G}}\ti{\sigma }$ est irr\'eductible pour tout $\lambda\in\Bbb R$ et donc avec (iii) l'op\'erateur d'entrelacement $J_{\ol{\ti{P}}\vert \ti{P}}(\ti{\sigma}_{\lambda\alpha})$ est toujours r\'egulier et bijectif, car non nul (cf. {\bf 2.3}). Maintenant, si l'op\'erateur d'entrelacement $J_{\ol{\ti{P}}\vert \ti{P}}$  est r\'egulier en $\ti{\sigma}$, alors par le (iii) du th\'eor\`eme la fonction $\mu $ est non nulle en $\ti{\sigma}$. Par ce qui pr\'ec\'edait, il suffit de consid\'erer le cas o\`u $w\ti{M}=\ti{M}$ et $w\ti{\sigma }=\ti{\sigma }$. Par hypoth\`ese, $\mu(\ti{\sigma })\ne 0$. Alors, par les propri\'et\'es (ii) a) et b) du th\'eor\`eme, $i_{\ti{P}}^{\ti{G}}\ti{\sigma }_{\lambda\alpha }$ est irr\'eductible pour tout $\lambda\in\Bbb R$, $\lambda\ne 0$. Comme avant, on en d\'eduit que $J_{\ol{\ti{P}}\vert \ti{P}}(\ti{\sigma}_{\lambda\alpha})$ est bijectif pour $\lambda\in\Bbb R$, $\lambda\ne 0$. Quant au cas $\lambda =0$, la bijectivit\'e de $J_{\ol{\ti{P}}\vert \ti{P}}(\ti{\sigma})$ a \'et\'e prouv\'ee dans la derni\`ere partie du lemme {\bf 3.1}. \hfill{\fin 2}

\null{\bf Annexe A: degr\'e formel et r\'esidu}

\null L'objet de cette annexe est la preuve de la partie (ii) de la proposition {\bf 2.4} qui est dans le cas d'un groupe r\'eductif contenu dans \cite{H2} utilisant \cite{H1} (cf. \cite{H2, p.51}). Silberger d\'eduit ce r\'esultant dans \cite{Si1} directement de la formule de Plancherel, mais comme la g\'en\'eralisation de \cite{H1, H2} a un int\'er\^et suppl\'ementaire, nous suivons \cite{H2}, en se restreignant toutefois aux parties essentielles.

\null{\bf A.1} On garde dans la suite les notations du paragraphe {\bf 1}. En particulier, $\ti{\sigma }$ sera une repr\'esentation cuspidale unitaire de $\ti{M}$, et on notera $\ti{\o }$ l'ensemble des classes d'isomorphie des repr\'esentations de la forme $\ti{\sigma}\otimes\ti{\chi }$ qui est un quotient de $\X (\ti{M})$, ce qui permettra de lib\'erer le symb\^ole $\ti{\sigma }$. Le stabilisateur dans $\X (\ti{M})$ d'un \'el\'ement de $\ti{\o }$ pour l'action par torsion ne d\'epend pas du choix de cet \'el\'ement et sera not\'e $Stab(\ti{\o })$. Le symb\^ole $\deg$ d\'esignera le degr\'e formel (voir par exemple \cite{W, III.1} et \cite{Li1, 2.5}) qui est d\'efini pour les repr\'esentations de carr\'e int\'egrable et est invariant par torsion par un caract\`ere non-ramifi\'e unitaire. Ceci permet d'\'etendre cette d\'efinition aux repr\'entations essentiellement de carr\'e int\'egrable.

\null
{\bf A.2} Si $\ti{\sigma }$ est dans $\ti{\o }$, on peut toujours trouver un caract\`ere non ramifi\'e $\ti{\chi}$ de $\ti{M}$ tel que $\ti{\sigma }\otimes\ti{\chi}^{-1}$ soit unitaire et \'ecrire $\ti{\chi }$ sous la forme $\ti{\chi}_{-\lambda\alpha}\chi_{\ti{G}}$ avec $\lambda\in\Bbb C$ et $\chi_{\ti{G}}$ dans $\X(\ti{G})$. La partie r\'eelle de $\lambda $ ne d\'epend pas du choix de $\ti{\chi }_{\ti{G}}$. On dira ici que $\Re(\ti{\sigma})=r$, si $\Re(\lambda)=r$ et $\chi_{\ti{G}}$ est unitaire.

On note $\Xim(\ti{M})$ (resp. $\Xim(A_M)$) le groupes des caract\`eres non ramifi\'es unitaires de $\ti{M}$ (resp. $A_M$). Le groupe $\Xim(A_M)$ est compact, et on le munit de la mesure de Haar de masse totale $1$. \'Ecrivons $\ti{\o}_0$ pour le sous-espace de $\ti{\o }$ form\'e des repr\'esentations unitaires. On munit $\Xim(\ti{M})$ (resp. $\ti{\o}_0$) de la mesure de Haar telle que la restriction $\Xim(\ti{M})\rightarrow \Xim(A_M)$ (resp. l'action de $\Xim(\ti{M})$ sur $\ti{\o}_0$) conserve localement les mesures. C'est conforme \`a \cite{H3, 1.6}, \cite{Li1, 2.1 et 2.6} et  \cite{W, I.1 et VI.3}. Si $r\in\Bbb R$ et $f$ est une fonction rationnelle sur $\ti{\o }$ telle que $\ti{\sigma}\mapsto f(\ti{\sigma}\otimes\ti{\chi}_{r\alpha})$ soit int\'egrable sur $\ti{\o} _0$, on \'ecrira $\int_{\Re(\ti{\sigma})=r}f(\ti{\sigma}) d\Im(\ti{\sigma})$ pour l'int\'egrale $\int_{\ti{\so }_0}   f(\ti{\sigma}\otimes\ti{\chi}_{r\alpha}) d\ti{\sigma }$.

Notons par ailleurs, pour $\ti{\sigma}\in\ti{\o }$, par $E_{\ti{P},\ti{\sigma }}^{\ti{G}}$ l'application lin\'eaire qui associe \`a un \'el\'ement $v\otimes v^{\vee }$ de $i_{\ti{P}\cap\ti{K}}^{\ti{K}}\ti{E}\otimes i_{\ti{P}\cap\ti{K}}^{\ti{K}}\ti{E}^{\vee }$ le coefficient matriciel $g\mapsto\langle (i_{\ti{P}}^{\ti{G}}\ti{\sigma})(g)v,v^{\vee }\rangle$ de la repr\'esentation induite. La notion d'une fonction lisse sur $\ti{\o }$ \`a valeurs dans $i_{\ti{P}\cap\ti{K}}^{\ti{K}}\ti{E}\otimes i_{\ti{P}\cap\ti{K}}^{\ti{K}}\ti{E}^{\vee }$ a \'et\'e explicit\'ee dans \cite{W, VI.1} et \cite{Li1, 2.6} (on demande l'invariance par l'action de $Stab(\ti{\o})$ qui y est pr\'ecis\'ee). On en d\'eduit instantan\'ement celle d'une fonction rationnelle ou polyn\^omiale.

\null{\bf Th\'eor\`eme:} \it Soit $\xi: \ti{\o }\rightarrow i_{\ol{\ti{P}}\cap\ti{K}}^{\ti{K}}\ti{E}\otimes i_{\ti{P}\cap\ti{K}}^{\ti{K}}\ti{E}^{\vee }$ une application polyn\^omiale \`a image dans un espace de dimension finie. Alors, la fonction $f_\xi:\ti{G}\rightarrow\Bbb C$ d\'efinie par
$$f_\xi(g)=\gamma(\ti{G}\vert\ti{M})^{-1}\int_{\Re(\ti{\sigma })=r\gg 0} E_{\ti{P},\ti{\sigma }}^{\ti{G}}((J_{\ti{P}\vert\ol{\ti{P}}}(\ti{\sigma})\otimes 1)\xi )(g^{-1})\mu(\ti{\sigma})\deg(\ti{\sigma}) d\Im(\ti{\sigma })$$ appartient \`a $C_c^{\infty}(\ti{G})$.

Posons $\varphi_\xi(\ti{\sigma})=(J_{\ti{P}\vert\ol{\ti{P}}}(\ti{\sigma })\otimes 1)\xi(\ti{\sigma })+(\lambda(\ti{w})\otimes J_{\ti{P}\vert\ol{\ti{P}}}(\ti{\sigma }^\vee )\lambda(\ti{w}))\xi(w^{-1}\ti{\sigma })$ et identifions $i_{\ti{P}\cap\ti{K}}^{\ti{K}}\ti{E}\otimes i_{\ti{P}\cap\ti{K}}^{\ti{K}}\ti{E}^{\vee }$ au sous-espace des applications de rang fini de $\End_\Bbb C(i_{\ti{P}\cap\ti{K}}^{\ti{K}}\ti{E})$. Supposons $w\ti{\o}=\ti{\o}$. Alors, pour tout $\ti{\sigma }\in\ti{\o }$, $(i_{\ti{P}}^{\ti{G}}\ti{\sigma })(f_\xi)=\varphi_\xi(\ti{\sigma})$.

\null Preuve: \rm Pour les groupes r\'eductifs, la premi\`ere partie de la proposition est \cite{H1, 2.1}. La preuve se g\'en\'eralise aux extensions centrales avec les remarques suivantes: il faut remplacer $A_0$ par $A_0^\dagger$ qui est d'indice fini dans $A_0$ \cite{Li1, 2.1.1}. Posons comme dans  \cite{H1, 2.1} $$M_0^+=\{m\in M_0\vert \forall\alpha\in\Delta \langle\alpha ,H_0(m)\rangle\geq 0\},$$ $A_0^+=A_0\cap M_0^+$ et ${A_0^\dagger}^+=A_0^\dagger\cap A_0^+$. Comme $A_0^\dagger$ est d'indice fini dans $A_0$, il existe un compact $C_0$ de $G$ tel que $\ti{C_0}\ti{{A_0^\dagger}^+}=\ti{A_0^+}$.

Partant de la d\'ecomposition de Cartan, tout se ram\`ene alors comme dans le cas r\'eductif \`a l'\'etude de fonctions sur $\ti{{A_0^\dagger}^+}$ de la forme $$a\mapsto \int_{\Re(\ti{\sigma } )=r\gg 0} p(\ti{\sigma })\langle i_{\ti{P}}^{\ti{G}}(\ti{\sigma})(a)J_{\ti{P}\vert\ti{\ol{P}}}(\ti{\sigma})v,v^\vee\rangle \mu(\ti{\sigma }) d\Im(\ti{\sigma })$$ avec $p$ polyn\^omial. Pour $G$ semi-simple, on pourra poursuivre la preuve comme dans le cas r\'eductif pour prouver que c'est \`a support compact, utilisant la proposition {\bf 2.5} \`a la place de \cite{H1, 1.3.1} et le fait que $H_M(A_M^\dagger)=H_{M_0}(A_0^\dagger)\cap a_M$ (cf. \cite{Li1, 2.1.1}).

Dans le cas g\'en\'eral, on a semblable au cas d'un groupe r\'eductif une d\'ecomposi-tion $\ti{A_0^\dagger}=\ti{A_G^\dagger} (\ti{A_0^\dagger\cap G^{der}})\ti{C_0'}$ avec $\ti{C_0'}$ compact, et on d\'ecompose $\Xim(\ti{M})$ comme dans le cas r\'eductif \`a l'aide de $\Xim(\ti{G})$ et de $\Xim(\ti{M\cap G^{der}})$. On peut alors ramener l'\'etude de la fonction sur $\ti{{A_0^\dagger}^+}$ qui est d\'efinie par une int\'egrale sur $\Xim(\ti{M})$ \`a celle d'une fonction sur $\ti{A_G^\dagger}$ d\'efinie par une int\'egrale sur $\Xim(\ti{G})$ et celle d'une fonction sur $\ti{A_0^\dagger\cap G^{der}}$ d\'efinie par une int\'egrale sur $\Xim(\ti{M\cap G^{der}})$. Il r\'esulte alors de la th\'eorie de Fourier sur un tore que l'int\'egrale sur $\Xim(\ti{G})$ d\'efinit une fonction \`a support compact sur $\ti{A_G^\dagger}$, alors que l'int\'egrale sur $\Xim(\ti{M\cap G^{der}})$ se traite relativement \`a $\ti{A_0^\dagger\cap G^{der}}$ comme dans le cas semi-simple.

La deuxi\`eme partie de la proposition correspond pour les groupes r\'eductifs \`a \cite{H1, 2.3}. Sa preuve est bas\'ee sur les propri\'et\'es de composition d'op\'erateurs d'entrelacement (qui sont les m\^emes pour les extensions centrales) et des formules int\'egrales standard. Elle reste donc valable pour les extensions centrales.\hfill{\fin 2}

\null {\bf A.3 Remarque:} Le th\'eor\`eme et sa preuve se g\'en\'eralisent pour $\ti{P}$ un sous-groupe parabolique standard arbitraire.

\null {\bf A.4} Notons $\iota_M$ l'indice de l'image de $\X(\ti{M})$ dans $\X(A_M)$ par l'application de restriction et analogue pour $\iota _G$. \'Ecrivons $\ti{\alpha}$ pour le multiple de $\alpha $ par une constante positive tel que $\langle H_M(h_\alpha),\ti{\alpha }\rangle=1$.\footnote{Dans \cite{H2}, ce multiple est not\'e $\alpha^*$. Le symb\^ole $h_{\alpha}$ d\'ej\`a introduit avant le lemme {\bf 1.2} a la m\^eme signification dans \cite{H2} et \cite{H3}.}

\null {\bf Proposition:} \it Soit $\ti{\pi }$ une repr\'esentation irr\'eductible de carr\'e int\'egrable de $\ti{G}$ dont le support cuspidal contient un \'el\'ement $\ti{\sigma }$ de $\ti{\o }$, $\Re(\ti{\sigma})>0$. Alors, $$\deg(\ti{\pi})={\iota_M\over\iota_G}\log q\ \gamma(\ti{G}\vert\ti{M})^{-1}\deg(\ti{\sigma })Res_{z=0}\mu(\ti{\sigma}_{z\ti{\alpha }}).$$

\null Preuve: \rm Remarquons d'abord qu'il r\'esulte de {\bf 3.2} (qui n'utilise pas {\bf 2.6}) que l'\'enonc\'e est vide si $w\ti{\sigma}\not\simeq\ti{\sigma}$. On pourra  donc supposer dans la suite $w\ti{\o }\simeq\ti{\o }$.

On se fixe une application polynomiale  $\xi: \ti{\o }\rightarrow i_{\ol{\ti{P}}\cap\ti{K}}^{\ti{K}}\ti{E}\otimes i_{\ti{P}\cap\ti{K}}^{\ti{K}}\ti{E}^{\vee }$ \`a image dans un espace de dimension finie. On va effectuer un changement de contours \`a l'expression int\'egrale pour $f_\xi $ et comparer le r\'esultat avec la formule de Plancherel.

Remarquons d'abord que, par construction, l'op\'erateur d'entrelacement et donc en particulier la fonction $\mu $ sont constantes sur les $\X(\ti{G})$-orbites et que le produit $\mu J_{\ti{P}\vert\ol{\ti{P}}}$ est r\'egulier si $\Re(\ti{\sigma })=0$ \cite{W, V.2.3}, \cite{Li1}.

Le changement de contours effectu\'e \`a l'expression pour $f_\xi $ donne par le th\'eor\`eme des r\'esidus, comme expliqu\'e dans \cite{H2, 3.6} (il faut adapter les notations et constants - qui sont $>0$ - au cas pr\'esent),
$$\eqalign {\gamma(\ti{G}\vert\ti{M})f_{\xi}(g)=&\int_{\Re(\ti{\sigma })=0} E_{\ti{P},\ti{\sigma }}^{\ti{G}}((J_{\ti{P}\vert\ol{\ti{P}}}(\ti{\sigma})\otimes 1)\xi(\ti{\sigma })) (g^{-1})\deg(\ti{\sigma })\mu(\ti{\sigma}) d\Im(\ti{\sigma })\cr +{\iota_M\over\iota_G}\log q\sum _i\int_{\ti{\so }_i}&E_{\ti{P},\ti{\sigma }}^{\ti{G}}((J_{\ti{P}\vert\ol{\ti{P}}}(\ti{\sigma})\otimes 1)\xi(\ti{\sigma })) (g^{-1})\deg(\ti{\sigma }) Res_{z=0}\mu(\ti{\sigma}_{z\ti{\alpha }})\ d\Im(\ti{\sigma}),\cr}$$ o\`u la somme porte sur les hyperplans r\'esiduels $\ti{\o_i}$ de $\mu $ qui sont des $\Xim(\ti{G})$-orbites de repr\'esentations $\ti{\sigma }_i\in\ti{\o }$ avec  $\Re(\ti{\sigma }_i)>0$. (La mesure sur $\ti{\o }_i$ est d\'efinie \`a partir de celle sur $\Xim(\ti{G})$ \cite{H2, 1.6}, demandant que l'action de $\Xim(\ti{G})$ pr\'eserve localement les mesures.)

Comme la fonction $\mu $ a un p\^ole en $\ti{\sigma }$, si $\ti{\sigma}\in\ti{\o}_i$, $J_{\ti{P}\vert\ol{\ti{P}}}(\ti{\sigma })$ est non-bijectif, r\'egulier et non nul par {\bf 2.3} et l'unique sous-repr\'esentation irr\'eductible de $i_{\ti{P}}^{\ti{G}}\ti{\sigma }$ est de carr\'e int\'egrable par {\bf 2.4} (i). Par cons\'equence, $E_{\ti{P},\ti{\sigma }}^{\ti{G}}((J_{\ti{P}\vert\ol{\ti{P}}}(\ti{\sigma})\otimes 1)\xi(\ti{\sigma })) (g^{-1})$ est une combinaison lin\'eaire de coefficients matriciels de cette repr\'esentation de carr\'e int\'egrable.

Effectuons le changement de variable $\ti{\sigma }\mapsto w^{-1}\ti{\sigma }$ sur l'int\'egrale dans le premier terme qui laisse invariant le degr\'e formel et la fonction $\mu $ (cf. \cite{W, V.2.1} et \cite{Li1, 2.4.3}, on trouve - apr\`es un petit calcul qui utilise la propri\'et\'e d'adjonction des op\'erateurs d'entrelacement \cite{W, IV.1 (11)}, \cite{Li1, 2.4.2 1.}\footnote{la preuve de cette propri\'et\'e d'adjonction n'est pas d\'etaill\'ee dans \cite{W, IV.1 (11)}, mais c'est une simple v\'erification utilisant les propri\'et\'es du quotient d'une mesure de Haar.} - que celle-ci est \'egale \`a $$\int_{\Re(\ti{\sigma })=0} E_{\ti{P},\ti{\sigma }}^{\ti{G}}((\lambda(\ti{w})\otimes J_{\ol{\ti{P}}\vert\ti{P}}(\ti{\sigma }^\vee )\lambda(\ti{w}))\xi(w^{-1}\ti{\sigma}))(g^{-1})\deg(\ti{\sigma }) \mu(\ti{\sigma}) d\Im(\ti{\sigma }).$$
Gr\^ace \`a la deuxi\`eme partie du th\'eor\`eme {\bf A.1}, on peut alors remplacer le premier terme par $${1\over 2}\int_{\Re(\ti{\sigma })=0} E_{\ti{P},\ti{\sigma }}^{\ti{G}}(\varphi_\xi(\ti{\sigma }))(g^{-1})\deg(\ti{\sigma })\mu(\ti{\sigma}) d\Im(\ti{\sigma }).$$

Remarquons que $\ti{\pi} (f_\xi )=0$ pour toute repr\'esentation de carr\'e int\'egrable $\ti{\pi}$ de $\ti{G}$ dont le support cuspidal n'est pas conjugu\'e \`a un \'el\'ement de $\ti{\o}$ (cf. \cite{H1, 2.4} qui se g\'en\'eralise aux extensions centrales pour les m\^emes raisons que \cite{H1, 2.3} utilis\'e dans la preuve du th\'eor\`eme {\bf A.2}) et que toute repr\'esentation irr\'eductible $\ti{\pi }$ de carr\'e int\'egrable de $\ti{G}$ de support cuspidal conjugu\'e \`a un \'el\'ement de $\ti{\o}$ peut se plonger dans une unique repr\'esentation de la forme $i_{\ti{P}}^{\ti{G}}\ti{\sigma }$ avec $\Re(\ti{\sigma})>0$ et $\ti{\sigma }$ dans $\ti{\o }$. Ainsi, $\ti{\pi}(f_\xi)$ est la restriction de $\varphi_\xi(\ti{\sigma })$ \`a l'espace de $\ti{\pi }$ que l'on notera dans la suite $\varphi_\xi (\ti{\pi})$.

Avec ces notations, la formule de Plancherel \cite{W, VIII.1.1}, \cite{Li1, 2.6} donne que $$\eqalign{f_\xi(g)=&\gamma(\ti{G}\vert\ti{M})^{-1}{1\over 2}\int_{\Re(\ti{\sigma })=0} E_{\ti{P},\ti{\sigma }}^{\ti{G}}(\varphi_\xi(\ti{\sigma }))(g^{-1})\deg(\ti{\sigma }) \mu(\ti{\sigma}) d\Im(\ti{\sigma })\cr &+\sum_i\deg(\ti{\pi}_i)\int_{\so_{\ti{\pi }_i}}E_{\ti{G}}^{\ti{G}}(\varphi_\xi (\ti{\pi}))(g^{-1})d\Im(\ti{\pi } ),\cr}$$ la somme portant sur les $\Xim(\ti{G})$-orbites $\o_{\ti{\pi }_i}$ de repr\'esentations de carr\'e int\'egrable $\ti{\pi }_i$ de $\ti{G}$ de support cuspidal dans $\ti{\o}$, les mesures \'etant celles pour laquelle l'action de $\Xim(\ti{G})$ pr\'eserve localement les mesures. Comparant les deux expressions pour $f_\xi$, on voit que les premiers termes - qui sont des int\'egrales sur $\ti{\o }_0$ - sont \'egaux.

Les autres expressions dans le premier terme sont des int\'egrales sur $\Xim(\ti{G})$ de coefficients matriciels de repr\'esentations de carr\'e int\'egrable. Par \cite{W, VII.2.2} (qui reste bien valable pour les extensions centrales - c'est implicite dans \cite{Li1}), des int\'egrales relatives \`a des orbites distinctes sont orthogonales. Les int\'egrales sur des orbites $\ti{\o}_i$ et $\ti{\o}_{\pi _i}$ qui se correspondent doivent donc \^etre \'egales, et l'application qui envoie $\ti{\pi}\in\ti{\o}_{\pi _i}$ sur l'\'el\'ement $\ti{\sigma }$ de $\ti{\o}_i$ qui est dans son support cuspidal est bien d\'efinie et bijective par la r\'eciprocit\'e de Frobenius qui montre \'egalement que cette bijection pr\'eserve les mesures. Reste \`a voir que, si $\ti{\sigma }$ et $\ti{\pi }$ se correspondent ainsi, alors $$E_{\ti{P},\ti{\sigma }}^{\ti{G}}((J_{\ti{P}\vert\ol{\ti{P}}}(\ti{\sigma})\otimes 1)\xi(\ti{\sigma }))=E_{\ti{G}}^{\ti{G}}(\varphi_\xi (\ti{\pi})).$$

En effet, $\varphi_\xi(\ti{\sigma })$ est par le th\'eor\`eme {\bf A.2} la somme de deux termes qui sont toutes les deux r\'eguli\`eres par {\bf 2.3} puisque $\Re(\ti{\sigma })>0$. Par la th\'eorie du quotient de Langlands (ou {\bf 3.2} qui ne d\'epend pas de cette annexe), l'image de $(\lambda(\ti{w})\otimes J_{\ti{P}\vert\ol{\ti{P}}}(\ti{\sigma }^\vee )\lambda(\ti{w}))\xi(w^{-1}\ti{\sigma })$ dans $End_{\Bbb C} (i_{\ti{P}\cap\ti{K}}^{\ti{K}}\ti{E})$ s'annule sur l'unique sous-repr\'esenta-tion irr\'eductible de $i_{\ti{P}}^{\ti{G}}\ti{\sigma }$ dont on sait qu'elle est isomorphe \`a $\ti{\pi }$. Ceci prouve l'\'egalit\'e ci-dessus.

Reste \`a remarquer que l'on peut choisir $\xi $ tel que les deux expressions soient non nulles: pour cela, il suffira par exemple de prendre pour $\xi $ la somme sur l'orbite par $Stab(\ti{\o})$ d'un \'el\'ement de la forme $v\otimes v^\vee $ avec $v$ dans l'espace de $\ti{\pi}$ et $v^\vee $ non nul sur $v$. On peut donc conclure que les constantes dans les deux int\'egrales sont \'egales, ce qui conduit \`a l'\'egalit\'e $$\deg(\ti{\pi})={\iota_M\over\iota_G}\log q\ \gamma(\ti{G}\vert\ti{M})^{-1}\deg(\ti{\sigma })Res_{z=0}\mu(\ti{\sigma}_{z\ti{\alpha }}).$$
\hfill{\fin 2}

\null{\bf Annexe B: op\'erateur autoadjoint}

\null On reprend les notations et hypoth\`eses du d\'ebut de la section {\bf 2.6}. En particulier, $\ti{\sigma }$ sera une repr\'esentation cuspidale unitaire de $\ti{M}$, v\'erifiant $\ti{w}\ti{\sigma}\simeq\ti{\sigma}$, et on fixe un produit scalaire $\ti{M}$ invariant $\langle\cdot,\cdot\rangle_{\ti{\sigma }}$ sur l'espace de $\ti{\sigma }$. On le supposera semi-lin\'eaire dans la premi\`ere composante et lin\'eaire dans la deuxi\`eme composante. L'objet de cette annexe est de d\'efinir l'isomorphisme $\rho_{\ti{\sigma },\ti{w}}:\ti{w}\ti{\sigma }\rightarrow \ti{\sigma }$ utilis\'e dans {\bf 2.6} dans la d\'efinition de l'op\'erateur d'entrelacement r\'egularis\'e $R_{\ti{P}}(\ti{\sigma}_{\lambda\alpha },\ti{w})$ dont nous montrerons qu'il est autoadjoint. Remarquons que l'on utilisera la proposition {\bf 3.1} qui ne d\'epend pas de {\bf 2.6}.

\null{\bf B.1} Si $\ti{\sigma }$ est un p\^ole de l'op\'erateur d'entrelacement standard $J_{\ol{\ti{P}}\vert\ti{P}}(\cdot)$, alors $\rho_{\ti{\sigma },\ti{w}}$ co\"\i ncid\'era avec l'op\'erateur avec la m\^eme notation d\'efini dans \cite{H3, 2.4}. Plus pr\'ecis\'ement, remarquant que la repr\'esentation induite $i_{\ti{P}}^{\ti{G}}\ti{\sigma }$ est irr\'eductible, $\rho_{\ti{\sigma },\ti{w}}$ sera l'unique isomorphisme qui induit par fonctorialit\'e l'isomorphisme $i_{\ti{P}}^{\ti{G}}w\ti{\sigma }\rightarrow i_{\ti{P}}^{\ti{G}}\ti{\sigma }$ donn\'e par $$[((Y(\ti{\chi })-1)\lambda(\ti{w})J_{\ol{\ti{P}}\vert \ti{P}}(\ti{\sigma}\otimes\ti{\chi }))_{\vert\ti{\chi}=1}]^{-1}$$ qui est bien d\'efini par {\bf 3.1}. Par abus de notation, on d\'esignera l'isomorphisme induit par fonctorialit\'e par le m\^eme symb\^ole. Observons que cet isomorphisme est \'egal \`a $\gamma(\ti{G}\vert\ti{M})^{-2}((Y(\ti{\chi })-1)^{-1}\mu(\ti{\sigma}\otimes\ti{\chi })J_{\ti{P}\vert \ol{\ti{P}}}(\ti{\sigma}\otimes\ti{\chi })\lambda(\ti{w}^{-1}))_{\vert\ti{\chi}=1}$ suite aux r\`egles de composition pour les op\'erateurs d'entrelacement et la d\'efinition de la fonction $\mu $ d\'ej\`a rappel\'ees dans les sections {\bf 1.} et {\bf 2.}.

Si $J_{\ol{\ti{P}}\vert \ti{P}}(\cdot)$ est r\'egulier en $\ti{\sigma}$, alors $\rho_{\ti{\sigma },\ti{w}}$ ne sera pas n\'ecessairement \'egal \`a l'isomor-phisme d\'efini dans \cite{H3, 2.5}. Pour le d\'efinir, on va d'abord travailler avec l'isomor-phisme d\'efini dans \cite{H3, 2.5} et le noter ici $\eta_{\ti{\sigma },\ti{w}}$. Plus pr\'ecis\'ement, $\eta_{\ti{\sigma },\ti{w}}$ sera un isomorphisme $w\ti{\sigma }\rightarrow \ti{\sigma }$ qui v\'erifie $\eta_{\ti{\sigma },\ti{w}}^2=\ti{\sigma}(\ti{w}^2)$, rappelant que $\ti{w}^2\in\ti{M}$. Ceci d\'etermine $\eta_{\ti{\sigma },\ti{w}}$ \`a un facteur $\pm 1$ pr\`es.

\null{\bf Lemme:} \it (i) Si $\ti{\sigma }$ est un p\^ole, l'adjoint de  $\rho_{\ti{\sigma },\ti{w}}$ relatif au produit scalaire $\langle\cdot,\cdot\rangle_{\ti{\sigma }}$ est $\rho_{w\ti{\sigma },\ti{w}^{-1}}$.

(ii) Si $\ti{\sigma }$ est un point r\'egulier, l'adjoint de $\eta_{\ti{\sigma },\ti{w}}$ relatif au produit scalaire $\langle\cdot,\cdot\rangle_{\ti{\sigma }}$ est \'egal \`a $\pm\eta_{\ti{\sigma },\ti{w}}^{-1}$.

\null Preuve: \rm (i) Soient $v,v'\in i_{\ti{P}\cap\ti{K}}^{\ti{K}}\ti{E}$. On pourra se restreindre ci-apr\`es \`a $\ti{\chi }$ \`a valeurs r\'eelles. Alors, utilisant que la contragr\'ediente de $\ti{\sigma }$ est isomorphe \`a la conjugu\'ee de $\ti{\sigma }$, l'invariance de la fonction $\mu $ par passage au dual et par $w$, et la propri\'et\'e d'adjonction de l'op\'erateur d'entrelacement (rappel\'ee dans la preuve de {\bf A.4}), on trouve $$\eqalign{\gamma(\ti{G}\vert\ti{M})^2\langle\rho_{\ti{\sigma },\ti{w}}v,v'\rangle_{\ti{\sigma }}&=\langle[(Y(\ti{\chi })-1)^{-1}\mu(\ti{\sigma}\otimes\ti{\chi })J_{\ti{P}\vert \ol{\ti{P}}}(\ti{\sigma}\otimes\ti{\chi })\lambda(\ti{w}^{-1})]_{\vert\ti{\chi}=1}v,v'\rangle_{\ti{\sigma }}\cr &=[(Y(\ti{\chi })-1)^{-1}\ol{\mu(\ti{\sigma}\otimes\ti{\chi })}\langle J_{\ti{P}\vert \ol{\ti{P}}}(\ti{\sigma}\otimes\ti{\chi })\lambda(\ti{w}^{-1})v,v'\rangle_{\ti{\sigma }}]_{\vert\ti{\chi}=1}\cr
&=[(Y(\ti{\chi })-1)^{-1}\mu(\ol{\ti{\sigma}}\otimes\ti{\chi }) \langle v,\lambda(\ti{w})J_{\ol{\ti{P}}\vert\ti{P}}(\ti{\sigma}\otimes\ti{\chi }^{-1})v'\rangle_{\ti{\sigma }}]_{\vert\ti{\chi}=1}\cr
&=[(Y(\ti{\chi })-1)^{-1}\mu(\ti{\sigma}\otimes\ti{\chi }^{-1}) \langle v,J_{\ti{P}\vert \ol{\ti{P}}}(\ti{w}\ti{\sigma}\otimes\ti{\chi })\lambda(\ti{w})v'\rangle_{\ti{\sigma }}]_{\vert\ti{\chi}=1}\cr
&=\gamma(\ti{G}\vert\ti{M})^2\langle v,\rho_{\ti{w}\ti{\sigma },\ti{w}^{-1}}v'\rangle_{\ti{\sigma }}.\cr
}$$

(ii) L'adjoint $\eta_{\ti{\sigma },\ti{w}}^*$ de $\eta_{\ti{\sigma },\ti{w}}$ est un isomorphisme $\ti{E}\rightarrow w\ti{E}$, comme on le v\'erifie facilement. Les  deux espaces \'etant irr\'eductibles, il existe un nombre complexe non nul $c$ tel que $\eta_{\ti{\sigma },\ti{w}}^*=c\eta_{\ti{\sigma },\ti{w}}^{-1}$. Or, par d\'efinition, $\eta_{\ti{\sigma },\ti{w}}^{-2}=\ti{\sigma}(\ti{w}^{-2})$ et $(\eta_{\ti{\sigma },\ti{w}}^*)^2=(\eta_{\ti{\sigma },\ti{w}}^2)^*$. Comme $\ti{w}^2\in\ti{M}$, l'adjoint de $\ti{\sigma}(\ti{w}^2)$ est $\ti{\sigma}(\ti{w}^{-2})$ par invariance du produit scalaire et donc $\ti{\sigma}(\ti{w}^{-2})=(\eta_{\ti{\sigma },\ti{w}}^*)^2=c^2\ti{\sigma}(\ti{w}^{-2})$. On en d\'eduit que $c=\pm 1$.\hfill{\fin 2}

\null{\bf B.2} Dans le cas o\`u $\ti{\sigma }$ est un point r\'egulier de $J_{\ol{\ti{P}}\vert\ti{P}}(\cdot)$, on pose $\rho_{\ti{\sigma },\ti{w}}=\eta_{\ti{\sigma },\ti{w}}$ si la constante $c$ dans le lemme ci-dessus vaut $1$ et $\rho_{\ti{\sigma },\ti{w}}=i\ \eta_{\ti{\sigma },\ti{w}}$ sinon.

Rappelons le produit scalaire $\langle\cdot, \cdot\rangle_{\ti{P},\ti{\sigma }}$ sur l'espace vectoriel $i_{\ti{P}\cap\ti{K}}^{\ti{K}}\ti{E}$ d\'eduit de $\langle\cdot,\cdot\rangle_{\ti{\sigma }}$ (cf. {\bf 2.6}).

\null {\bf Proposition:} \it L'op\'erateur d'entrelacement $R_{\ti{P}}(\ti{\sigma}_{\lambda\alpha },\ti{w})$ est, pour tout $\lambda\in\Bbb R$, autoadjoint relatif \`a $\langle \cdot,\cdot\rangle_{\ti{P},\ti{\sigma }}$.

\null Preuve: \rm On va distinguer les cas $\ti{\sigma }$ p\^ole et point r\'egulier de $J_{\ol{\ti{P}}\vert\ti{P}}(\cdot)$.

(i) Supposons que $\ti{\sigma }$ soit un p\^ole. Par prolongation analytique, il suffira de consid\'erer le cas $\lambda\ne 0$. Dans ce cas, le facteur $X-1$ peut \^etre omis. Utilisant les propri\'et\'es d'adjoints prouv\'ees et remarqu\'ees dans {\bf B.1} et les propri\'et\'es basiques de $\rho_{\ti{\sigma },\ti{w}}$ d\'ej\`a utilis\'e dans \cite{H3, 2.}, on trouve, avec $v,v'\in i_{\ti{P}\cap\ti{K}}^{\ti{K}}\ti{E}$,
$$\eqalign{\langle \rho_{\ti{\sigma },\ti{w}}\lambda(\ti{w})J_{\ol{\ti{P}}\vert\ti{P}}(\ti{\sigma}_{\lambda\alpha })v,v'\rangle_{\ti{P},\ti{\sigma }}&=\langle v,J_{\ti{P}\vert\ol{\ti{P}}}(\ti{\sigma}_{-\lambda\alpha })\lambda(\ti{w}^{-1})\rho_{\ti{w}\ti{\sigma },\ti{w}^{-1}}v'\rangle_{\ti{P},\ti{\sigma }}\cr &=\langle v,\lambda(\ti{w}^{-1})J_{\ol{\ti{P}}\vert\ti{P}}(\ti{w}\ti{\sigma}_{\lambda\alpha })\rho_{\ti{w}\ti{\sigma },\ti{w}^{-1}}v'\rangle_{\ti{P},\ti{\sigma }}\cr
&=\langle v,\lambda(\ti{w}^{-1})\rho_{\ti{w}\ti{\sigma },\ti{w}^{-1}}J_{\ol{\ti{P}}\vert\ti{P}}(\ti{\sigma}_{\lambda\alpha })v'\rangle_{\ti{P},\ti{\sigma }}.\cr
}$$
On est donc ramen\'e \`a montrer que $\lambda(\ti{w}^{-1})\rho_{\ti{w}\ti{\sigma },\ti{w}^{-1}}=\rho_{\ti{\sigma },\ti{w}}\lambda(\ti{w})$. Pour cela, on pourra montrer que leurs compos\'es \`a droite avec $\rho_{\ti{w}\ti{\sigma },\ti{w}^{-1}}$ sont \'egaux. Par proposition \cite{H3, 2.4},
$$\gamma(\ti{G}\vert\ti{M})^2\rho_{\ti{w}\ti{\sigma },\ti{w}^{-1}}^2=[(Y(\ti{\chi })-1)^{-1} (Y(\ti{\chi })^{-1}-1)^{-1}\mu(\ti{\sigma}\otimes\ti{\chi })]_{\ti{\chi }=1}\lambda(\ti{w}^2).$$
Remarquons que $\lambda(\ti{w}^{-1})\rho_{\ti{w}\ti{\sigma },\ti{w}^{-1}}\lambda(\ti{w})=\rho_{\ti{w}\ti{\sigma },\ti{w}^{-1}}$. Avec cela, on trouve
$$\eqalign{&\gamma(\ti{G}\vert\ti{M})^4\rho_{\ti{\sigma },\ti{w}}\lambda(\ti{w})\rho_{\ti{w}\ti{\sigma },\ti{w}^{-1}}\cr =&[(Y(\ti{\chi })-1)^{-1}\mu(\ti{\sigma}\otimes\ti{\chi })J_{\ti{P}\vert \ol{\ti{P}}}(\ti{\sigma}\otimes\ti{\chi })\lambda(\ti{w}^{-1})\cr &\qquad\qquad\qquad(Y(\ti{\chi }^{-1})-1)^{-1}\mu(\ti{w}\ti{\sigma}\otimes\ti{\chi }^{-1})J_{\ti{P}\vert \ol{\ti{P}}}(\ti{w}\ti{\sigma}\otimes\ti{\chi }^{-1})\lambda(\ti{w})^2]_{\vert\ti{\chi}=1}\cr
=&[(Y(\ti{\chi })-1)^{-1}(Y(\ti{\chi }^{-1})-1)^{-1}\mu(\ti{\sigma}\otimes\ti{\chi })\mu(\ti{\sigma}\otimes\ti{\chi })\cr &\qquad\qquad\qquad\qquad\qquad\qquad\qquad\qquad J_{\ti{P}\vert \ol{\ti{P}}}(\ti{\sigma}\otimes\ti{\chi })J_{\ol{\ti{P}}\vert\ti{P}}(\ti{\sigma}\otimes\ti{\chi })\lambda(\ti{w})]_{\vert\ti{\chi}=1}\cr
=&\gamma(\ti{G}\vert\ti{M})^2[(Y(\ti{\chi })-1)^{-1}(Y(\ti{\chi }^{-1})-1)^{-1}\mu(\ti{\sigma}\otimes\ti{\chi })]_{\vert\ti{\chi}=1}\lambda(\ti{w})\cr
}.$$
On a donc bien \'egalit\'e, et c'est autoadjoint.

(ii) Supposons maintenant que $\ti{\sigma }$ soit un point r\'egulier. Analogue au cas pr\'ec\'e-dent, on trouve
$$\langle \rho_{\ti{\sigma },\ti{w}}\lambda(\ti{w})J_{\ol{\ti{P}}\vert\ti{P}}(\ti{\sigma}_{\lambda\alpha })v,v'\rangle_{\ti{P},\ti{\sigma }}=\langle v,\rho_{\ti{\sigma },\ti{w}}^*\lambda(\ti{w}^{-1})J_{\ol{\ti{P}}\vert\ti{P}}(\ti{\sigma}_{\lambda\alpha })v'\rangle_{\ti{P},\ti{\sigma }}.$$
Il faut donc ici montrer que $\rho_{\ti{\sigma },\ti{w}}^*\lambda(\ti{w}^{-1})=\rho_{\ti{\sigma },\ti{w}}\lambda(\ti{w})$, ce qui \'equivaut \`a $\lambda(\ti{w}^{-2})=(\rho_{\ti{\sigma },\ti{w}}^*)^{-1}\rho_{\ti{\sigma },\ti{w}}$. Remarquons que $\lambda(\ti{w}^{-2})$ est induit par l'automorphisme $\ti{\sigma }(\ti{w}^2)$ de l'espace vectoriel $\ti{E}$. Quant \`a $(\rho_{\ti{\sigma },\ti{w}}^*)^{-1}\rho_{\ti{\sigma },\ti{w}}$, il faut distinguer deux cas. Si $\rho_{\ti{\sigma },\ti{w}}^*=\eta_{\ti{\sigma },\ti{w}}^{-1}$, on trouve bien l'\'egalit\'e par le lemme {\bf B.1} et la d\'efinition de $\rho_{\ti{\sigma },\ti{w}}$. Dans l'autre cas, on observe que $\rho_{\ti{\sigma },\ti{w}}^*=-i \eta_{\ti{\sigma },\ti{w}}^*=i\eta_{\ti{\sigma },\ti{w}}^{-1}=-\rho_{\ti{\sigma },\ti{w}}^{-1}$, d'o\`u $(\rho_{\ti{\sigma },\ti{w}}^*)^{-1}\rho_{\ti{\sigma },\ti{w}}=-\rho_{\ti{\sigma },\ti{w}}^2=\eta_{\ti{\sigma },\ti{w}}^2=\ti{\sigma }(\ti{w}^2)$.\hfill{\fin 2}

\enddocument